\documentclass[a4paper,11pt]{article}

\usepackage{enumerate}

\usepackage[margin=2cm,footskip=30pt]{geometry}
\usepackage{layout}

\usepackage[utf8]{inputenc}
 \usepackage{multirow}
\usepackage[T1]{fontenc,url}
\usepackage[english]{babel}
\usepackage{tikz}
\usepackage{graphicx}
\usepackage[english]{babel}
\usepackage{setspace}
\usepackage{color}
\usepackage{amsmath,amsthm,amsfonts,mathrsfs}
\usepackage{csquotes,textcomp}
\usepackage{hyperref} 
\usepackage{times}
\usepackage{tikz}
\usetikzlibrary{trees}
\usetikzlibrary{arrows,}
\usepackage{stmaryrd}
\usepackage{graphicx}
\usepackage{subfigure,float}
\usepackage{placeins}

\usepackage{algorithm}
\usepackage[noend]{algorithmic}

\usepackage{amssymb}
\usepackage{mathptmx}
\usepackage{authblk}
\usepackage{bbm}
\usepackage{mathtools}

\usepackage[toc,page]{appendix}

\date{ \mbox{} }
\newcommand\myeq{\stackrel{\mathclap{\normalfont\mbox{D}}}{=}}

\begin{document}
    
\title{Parameter Estimation in Gaussian Mixture Models \\ with Malicious Noise, \\ without Balanced Mixing Coefficients}
\author{ Jing Xu, Jakub Marecek  }
\maketitle


    

\newtheorem{thm}{Theorem}[section]
\newtheorem{exa}{Example}[section]
\newtheorem{lemma}[thm]{Lemma}
\newtheorem{defi}{Definition}
\newtheorem{cor}[thm]{Corollary}
\newtheorem{prop}[thm]{Proposition}
\newcommand{\PR}{\mathbb{P}}
\newcommand{\E}{\mathbb{E}}

\newcommand{\tr}{\mathrm{tr}}
\newcommand{\p}{\mathbf{p}}
\newcommand{\q}{\mathbf{q}}
\newcommand{\bv}{\mathbf{v}}
\newcommand{\x}{\mathbf{x}}
\newcommand{\y}{\mathbf{y}}

\newcommand\righttwoarrow{%
        \mathrel{\vcenter{\mathsurround0pt
                \ialign{##\crcr
                        \noalign{\nointerlineskip}$\rightarrow$\crcr
                        \noalign{\nointerlineskip}$\rightarrow$\crcr
                }%
        }}%
}

\begin{abstract}
We consider the problem of estimating the means of two Gaussians in a 2-Gaussian mixture, which is not balanced and 
 is corrupted by noise of an arbitrary distribution. 
We present a robust algorithm to estimate the parameters,
together with upper bounds on the numbers of samples required for 
the estimate to be correct,
where the bounds are parametrised by the dimension, ratio of the mixing coefficients,
a measure of the separation of the two Gaussians, related to Mahalanobis distance,
and a condition number of the covariance matrix.
In theory, this is the first sample-complexity result for imbalanced mixtures corrupted by adversarial noise.
In practice, our algorithm outperforms the vanilla
Expectation-Maximisation (EM) algorithm
in terms of estimation error. 
\end{abstract}

\section{Introduction}

Gaussian mixture models are central to both theory and practice of Statistics \cite{titterington1985statistical,McLachlan2000}.
As a result of more than a century of study \cite{Pearson1894},
there are algorithms \cite{Hardt2015} in the noise-free setting with balanced mixing coefficients, 
which are essentially optimal \cite{Yang1999} with respect to their sample complexity and time complexity.
Even for the expectation-maximisation algorithm, which is often used in practice \cite{titterington1985statistical},
there is now some understanding of the performance \cite{Faria2010,Caramanis2014am,balakrishnan2017,caichime} and its limitations.
Within robust statistics \cite{Huber1964} and agnostic learning \cite{lai2016agnostic,Diakonikolas2017b},
there has been recent progress in estimating parameters of a single Gaussian from a mixture of the Gaussian and noise.
Within mixed regression, there has been some progress in estimating parameters of a mixture \cite{Chen2014,hand2016convex} with balanced coefficients,
but very little \cite{7045499,naim2012convergence} is known otherwise.

We propose to study the problem that relaxes both the noise-free assumption and the assumption on the balance of the mixing coefficients:

\begin{defi}[Robust Parameter Estimation in Noisy 2-GMM]  \label{def:P1}
Given $m$ points in $\mathbb{R}^n$ that are each, with probability $w_1 > 0$ from an unknown 
Gaussian $N(\mu_1, \Sigma)$, with probability $w_2 > 0$ from an unknown Gaussian distribution $N( \mu_2,\Sigma)$, and 
with probability $w_3 = 1-w_1 - w_2 > 0$ completely arbitrary, estimate $\mu_1, \mu_2$ and $\Sigma$.
\end{defi}
Throughout the paper, we assume:

\noindent \textbf{Assumption 1.} $w_1 > w_2 > w_3$.

Moreover,
\begin{itemize}
\item we consider arbitrary, adversarial noise. Our sample complexity is parametrised by the proportion $w_3$ of noise among the samples.
\item We do not make further assumptions on the balance between the Gaussians. Instead, our results are parametrised by 
the ratios of mixing coefficients $w_2/w_1$ and $w_3/w_1$.
\end{itemize}

\begin{figure}[t]
\centering
        \hskip -0.1cm
        \begin{minipage}{0.49\textwidth}
            \includegraphics[width=1\textwidth,trim={1cm 5cm 1cm 5cm},clip]{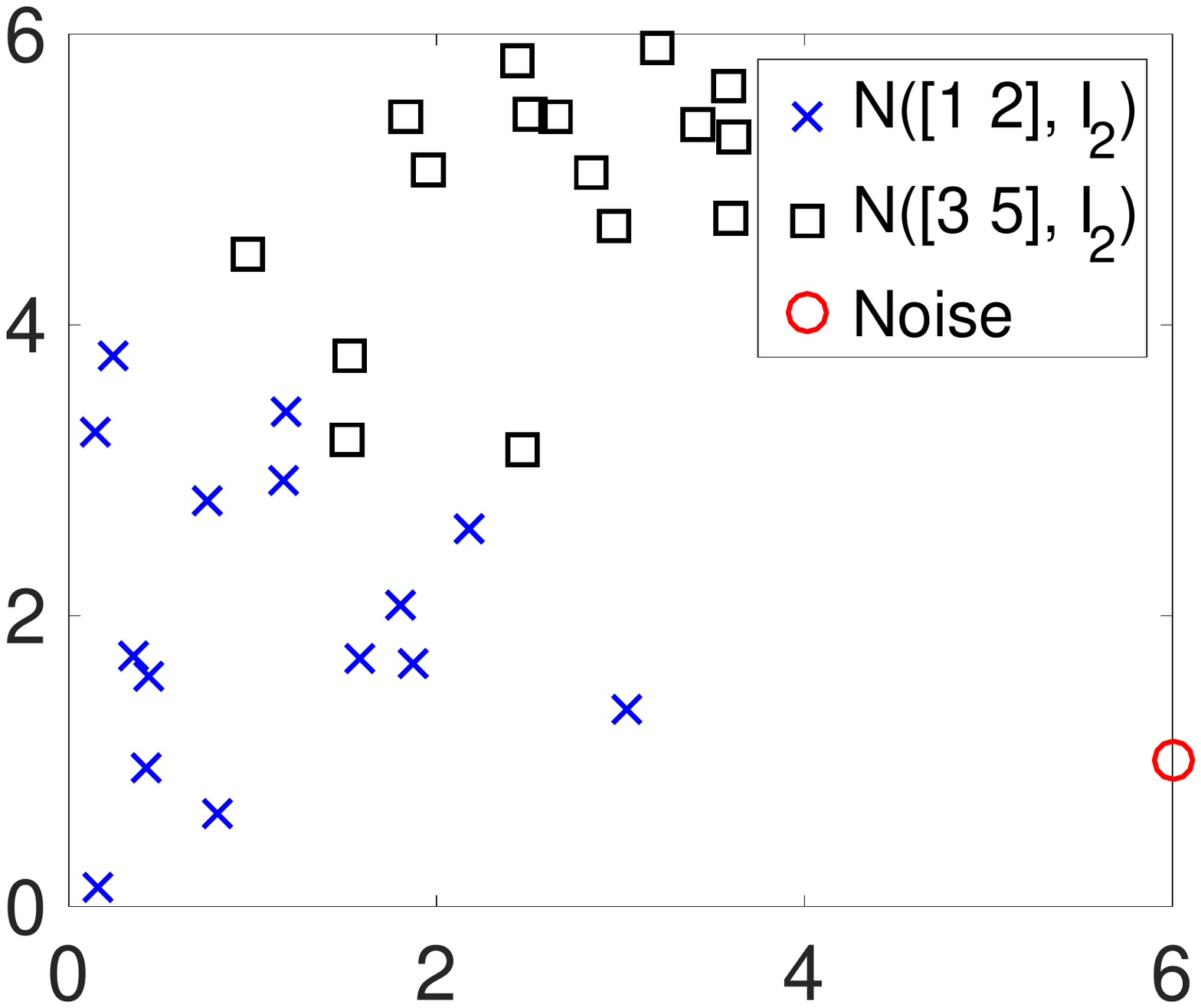} 
                \end{minipage}
        \hskip -0.1cm
        \begin{minipage}{0.49\textwidth}
            \includegraphics[width=1\textwidth,trim={1cm 5cm 1cm 5cm},clip]{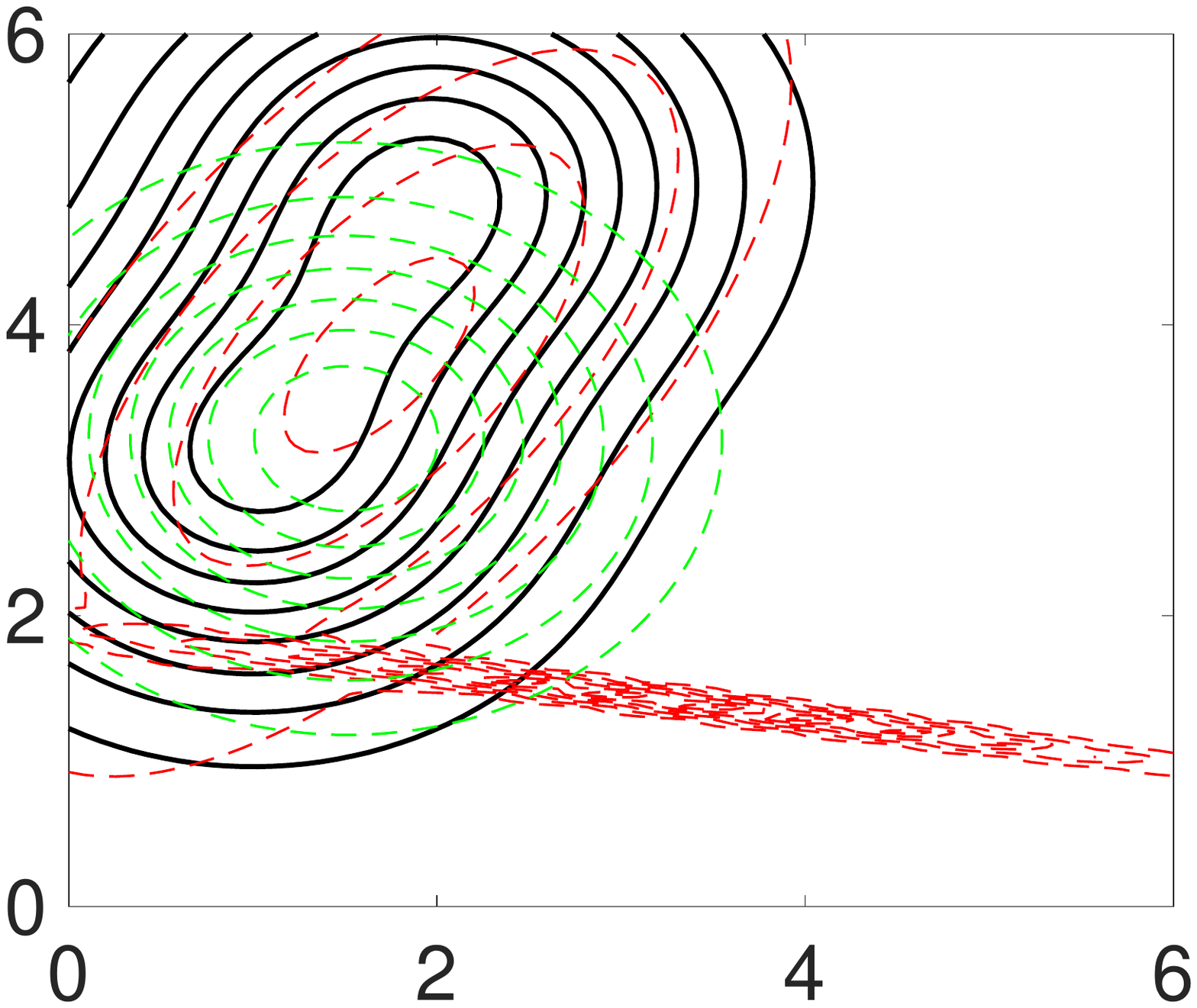} 
                \end{minipage}
                    \caption{Left: 20 and 20 samples from two Gaussians and one sample of noise at $(6,1)$.
                    Right: A contour plot of the 2-GMM (in black solid lines),
                    an EM estimate from the 41 samples (in red dashed lines),
                    and an estimate by Algorithm~\ref{alg1} using the same 41 samples (in green dashed lines).}
                    \label{motivation1}
        \end{figure}

\begin{figure}[t]
    \centering
    \hskip -0.1cm    
        \begin{minipage}{0.49\textwidth}
            \includegraphics[width=1\textwidth,trim={1cm 5cm 1cm 5cm},clip]{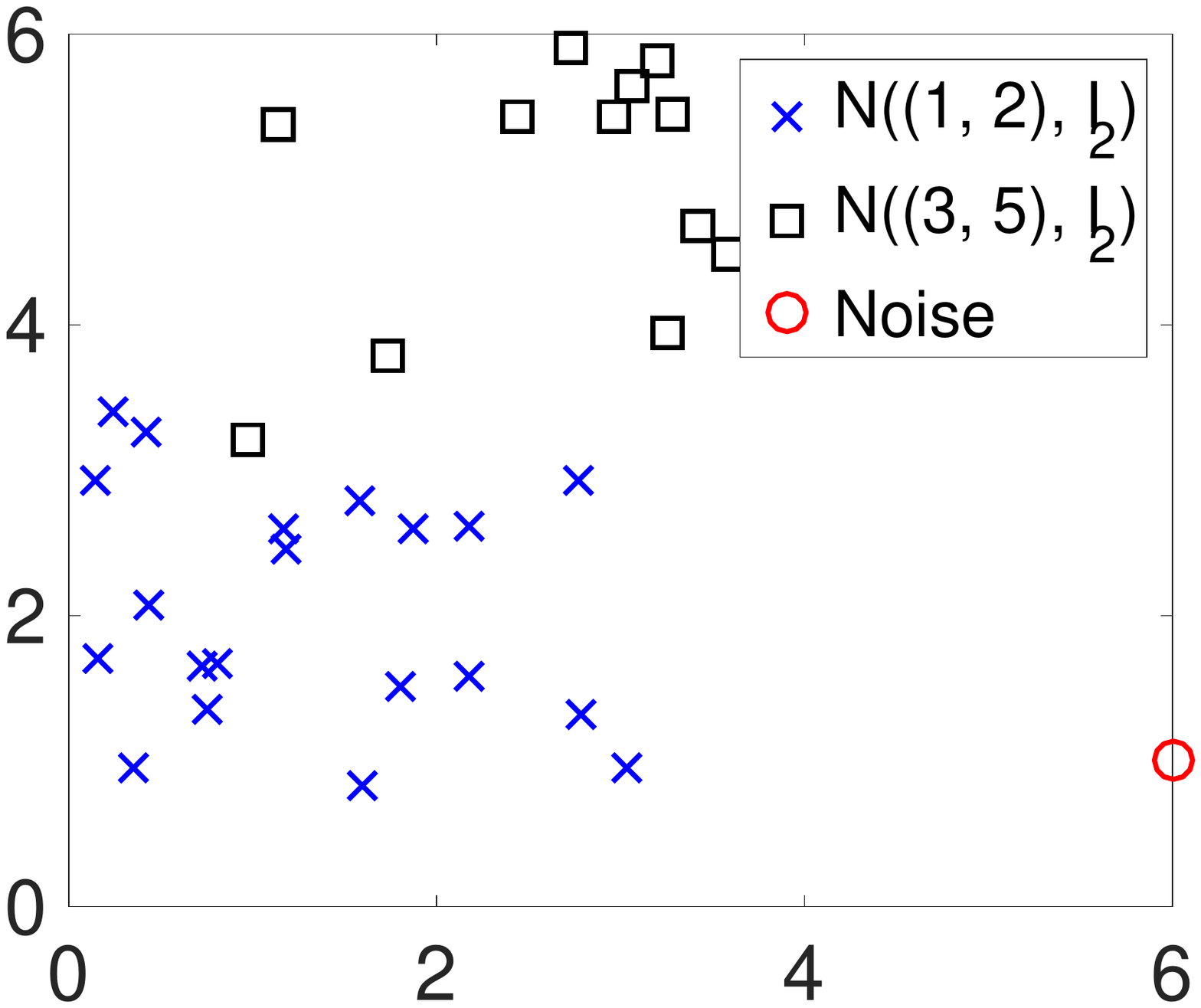}
                \end{minipage}
        \hskip -0.1cm
        \begin{minipage}{0.49\textwidth}
            \includegraphics[width=1\textwidth,trim={1cm 5cm 1cm 5cm},clip]{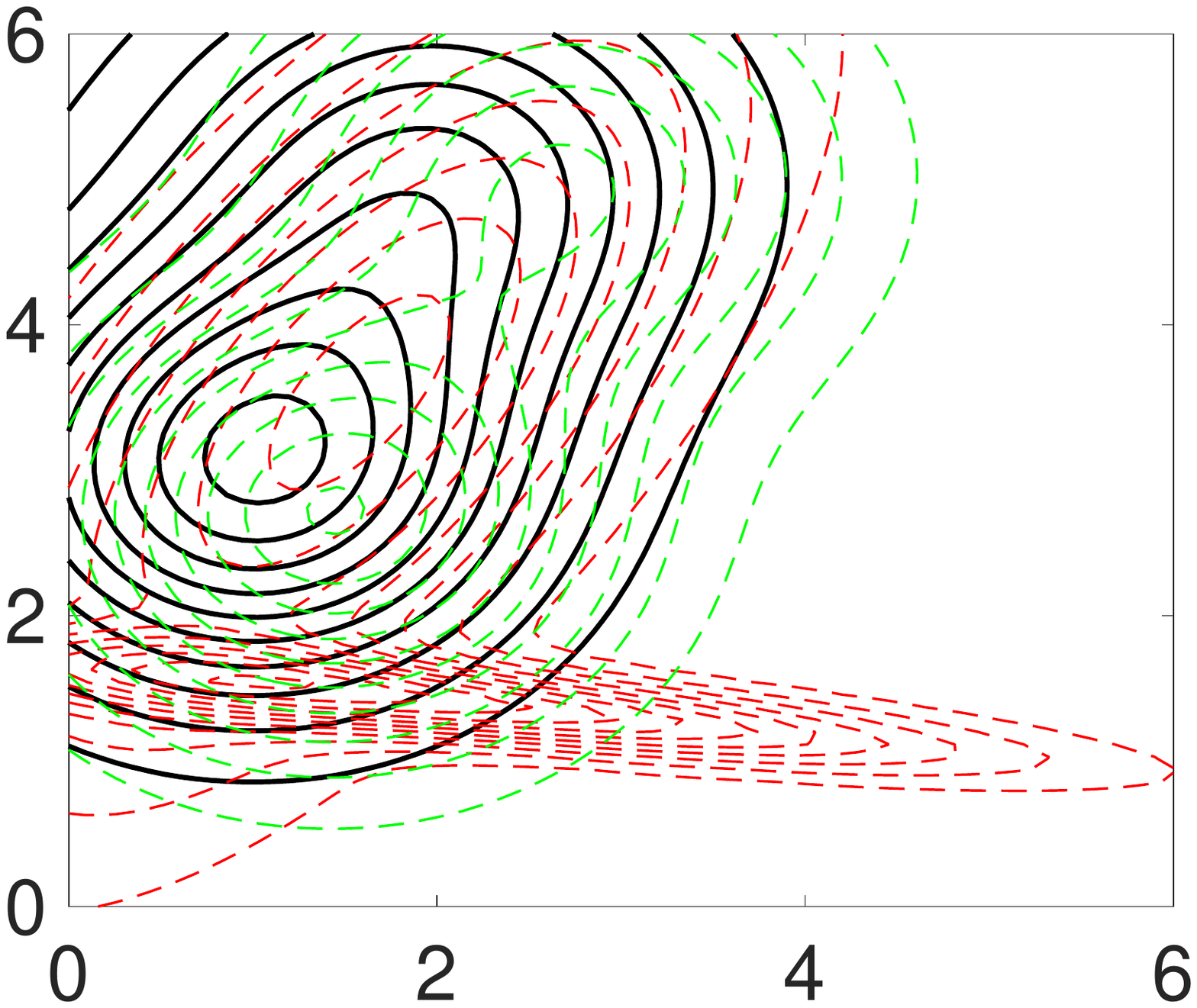}
                \end{minipage} \\
    \hskip 0.1cm    
        \begin{minipage}{0.49\textwidth}
            \includegraphics[width=1\textwidth,trim={1cm 5cm 1cm 5cm},clip]{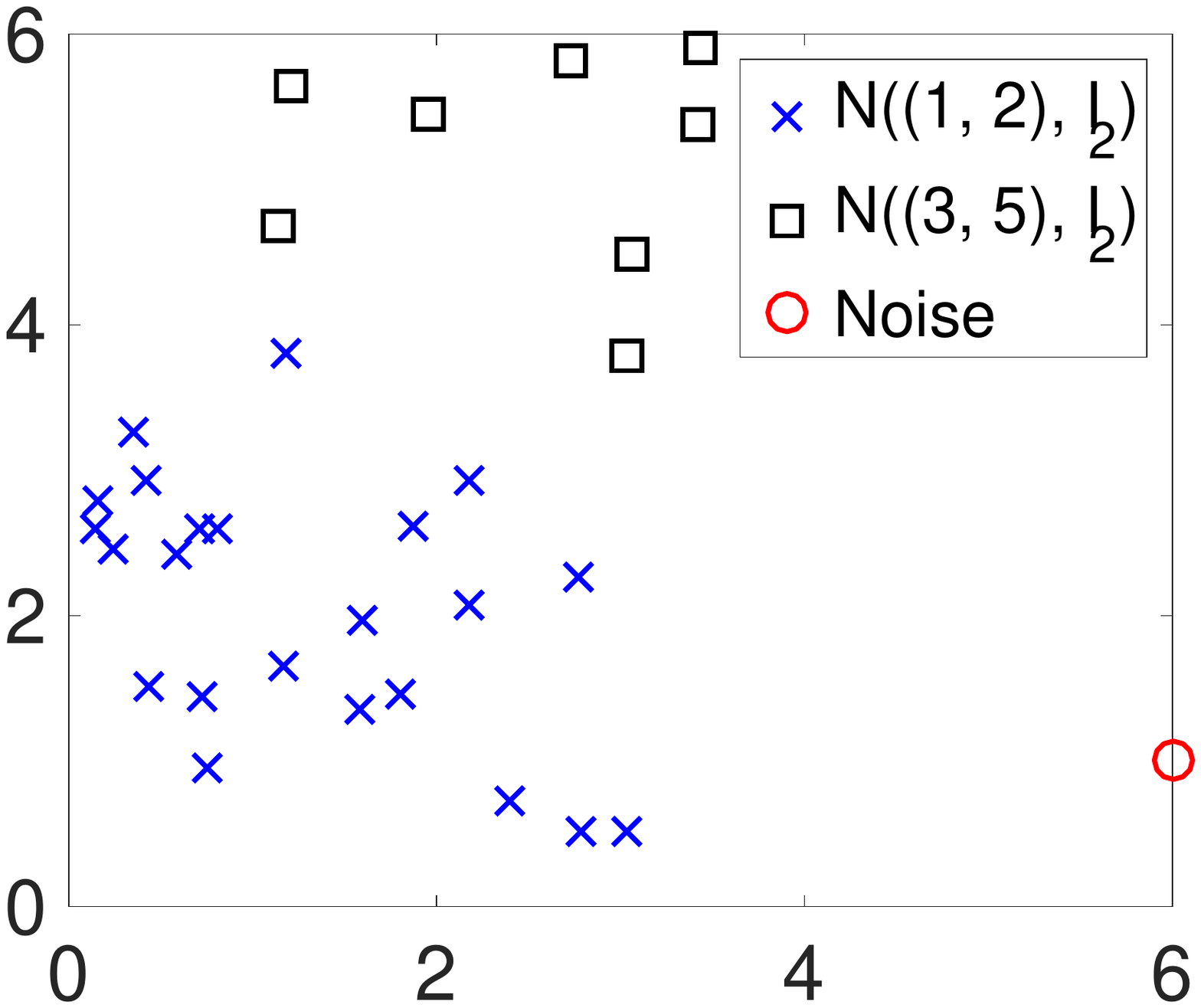}
                \end{minipage}
        \hskip -0.1cm
        \begin{minipage}{0.49\textwidth}
            \includegraphics[width=1\textwidth,trim={1cm 5cm 1cm 5cm},clip]{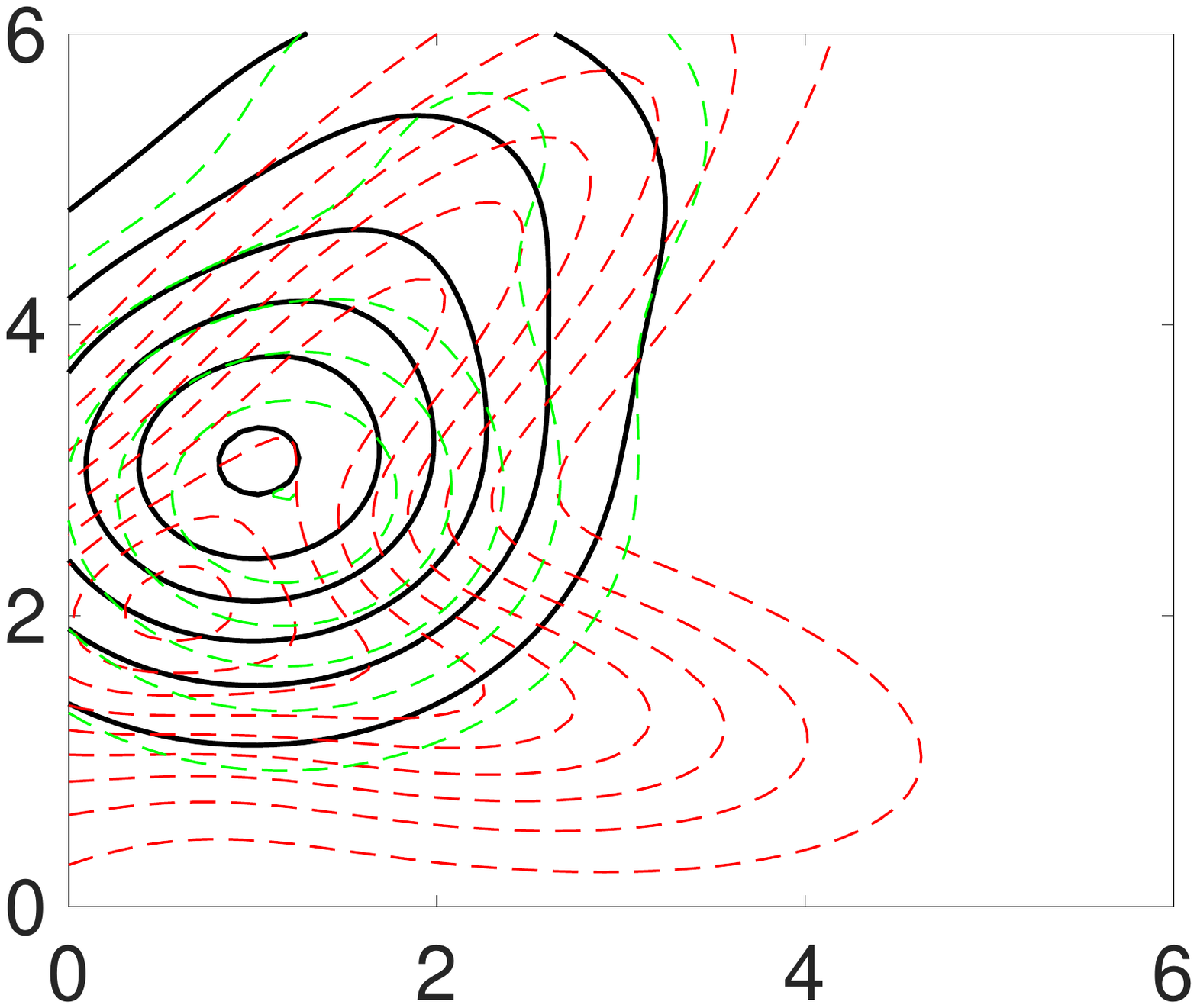}
                \end{minipage} \\                
    \hskip 0.1cm    
        \begin{minipage}{0.49\textwidth}
            \includegraphics[width=1\textwidth,trim={1cm 5cm 1cm 5cm},clip]{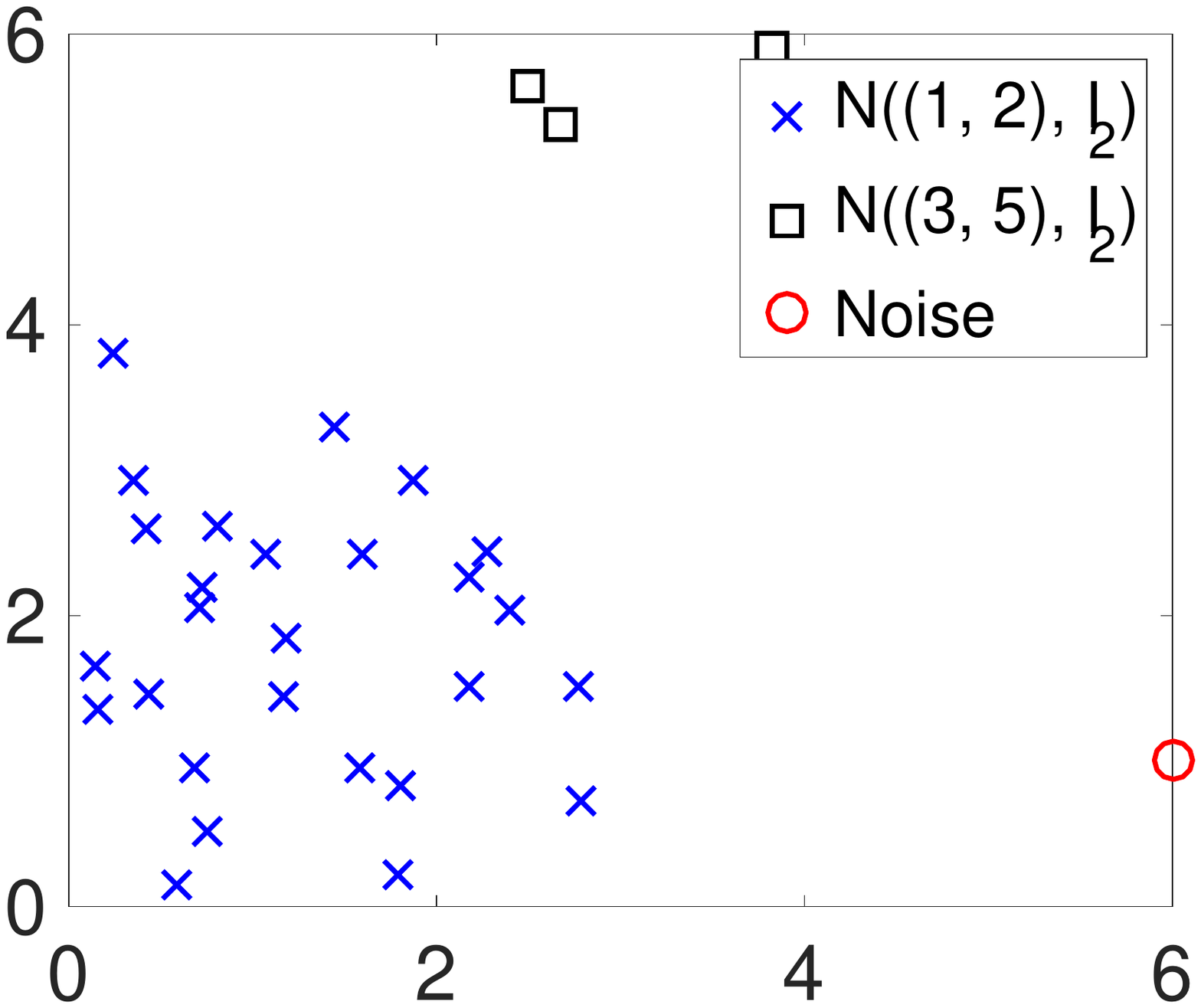}
                \end{minipage}
        \hskip -0.1cm
        \begin{minipage}{0.49\textwidth}
            \includegraphics[width=1\textwidth,trim={1cm 5cm 1cm 5cm},clip]{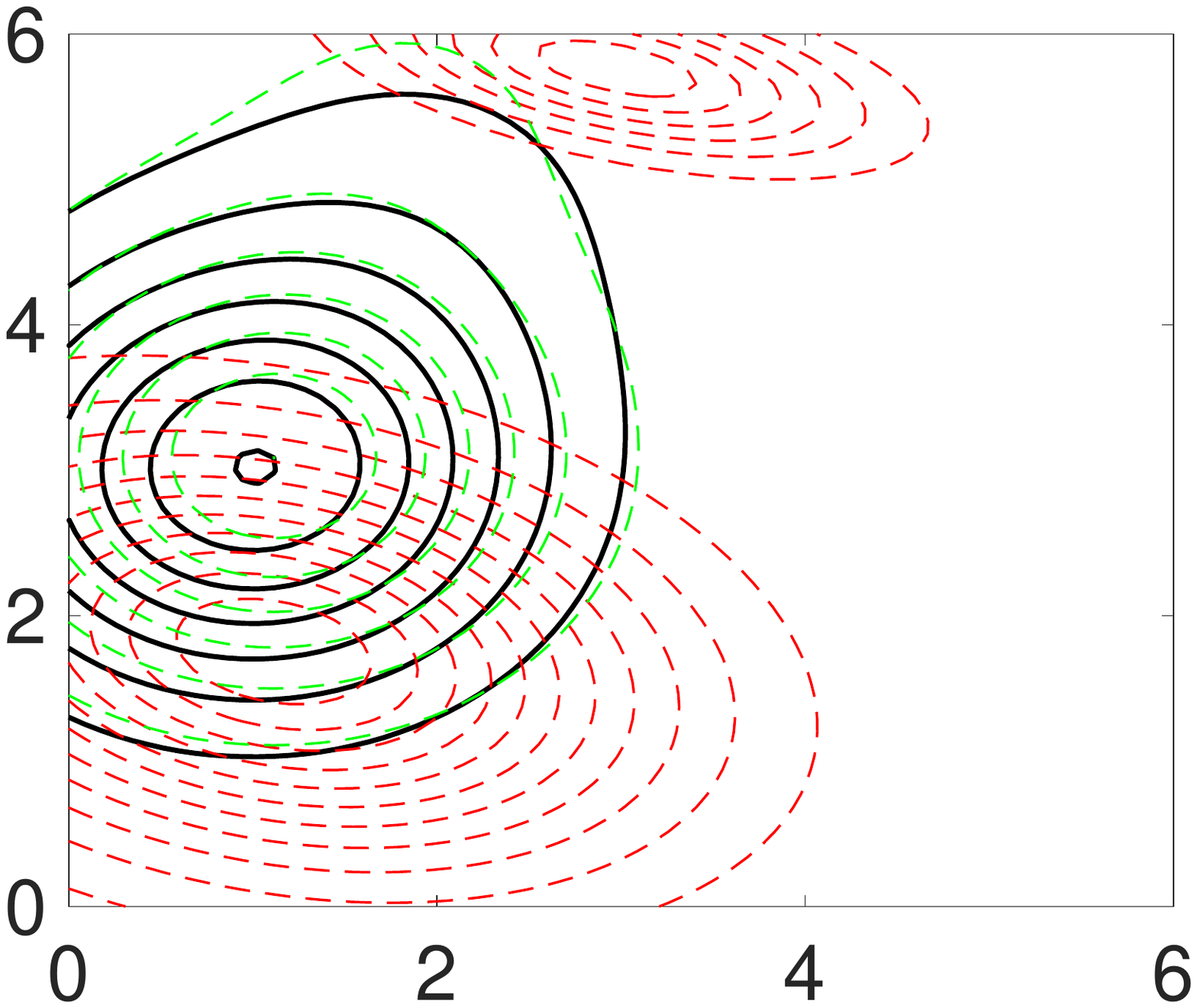}
                \end{minipage} \\                      
                    \caption{Left: $25$ (top), $30$ (middle), and $35$ (bottom) samples from one Gaussian and $15, 10$, and $5$ samples from the other Gaussian, 
                    with one sample of noise at $(6,1)$. 
                    Right: A contour plot of the 2-GMM (in black solid lines),
                    an EM estimate from the 41 samples (in red dashed lines),
                    and an estimate of Algorithm~\ref{alg1} using the same 41 samples (in green dashed lines).}
                    \label{motivation2}                
\end{figure}

The importance of \emph{not} making further assumptions is hard to overstate. 
For simplicity, let us illustrate this on a noisy mixture of $N((1, 2), I_2)$ and
$N((3, 5), I_2)$ in $\mathbb{R}^2$, with $I_2$ being the $2 \times 2$ identity matrix. 
First, we consider the situation, where the mixture is balanced, 
$m = 41, w_1 = 20/41, w_2 = 20/41$, and there is a single sample of noise at $(6,1)\in \mathbb{R}^2$.
Figure~\ref{motivation1} compares the estimate obtained using a standard
expectation-maximisation (EM) algorithm 
with the true mixture and our algorithm.
Notice that the one sample of noise causes the EM algorihm to mis-estimate the second component completely.
Figure~\ref{motivation2} illustrates the impact of varying mixing coefficients. 
In the three rows of sub-figures, we have used $w_1 = 25/41, 30/41$, and $35/41$, respectively, 
and the same single element of noise throughout, $w_3 = 1/41$. 
For $w_1 > 35/41$, the performance of the EM algorithm deteriorates quickly.
We stress that the samples have been obtained with the random seed set to 1 
and are not pathological.
The EM algorithm is implemented by \textit{fitgmdist} in MathWorks Matlab 2016a,
with default parameters.
One could conclude that algorithms designed for the noise-free, balanced case
should be avoided in applications, 
where the assumptions may be easily violated.

For robust parameter estimation in  noisy 2-GMM,
we present an iterative algorithm.
The algorithm could be seen as the extension of algorithms for 
the estimation of the mean and covariance of a single Gaussian in the presence of malicious noise 
(noisy single gaussian model),
as studied by \cite{lai2016agnostic}, 
to the noisy 2-GMM and beyond.
In our proposed algorithm, 
each iteration considers one Gaussian, in the decreasing order of their mixing coefficients.
In each iteration, parameters of one Gaussian are estimated under the assumption that the remaining
 samples are either from the Gaussian or 
the arbitrarily-distributed noise.
At the end of each iteration, the samples corresponding to the Gaussian are filtered out.
This way, Robust Parameter Estimation in Noisy 2-GMM can be reduced to 
2 calls of an algorithm for parameter estimation in Noisy 1-GMM and 
some additional processing in time $O(mn^2 + m \log m)$.


Outside of the meta-algorithm, our contributions are as follows:
\begin{itemize}
\item We prove an upper bound on the number of samples required to reach a given precision, 
considering a spectral method of \cite{lai2016agnostic}, with the bound parametrised by the dimension $n$,
 ratios of the mixing coefficients $w_2/w_1$ and $w_3/w_1$, a condition number of the covariance matrix, and a measure of the separation of the two Gaussians, related to Mahalanobis distance. 
\item In both theory and computational illustrations, we show that the algorithm is surprisingly robust to the error in the input $w_1 \in [0,1]$. 
\item In computational tests with Cauchy-distributed noise $w_3 > 0$, we demonstrate that the performance of the algorithm
 employing a spectral method of \cite{lai2016agnostic}
 is superior to the vanilla expectation-maximisation algorithm in terms of estimation error.
\end{itemize}

\section{The Algorithm}

We present an algorithm that estimates the parameters of Gaussians 
in the decreasing order of their mixing coefficients.
For each Gaussian component in turn, an algorithm for estimating the mean and covariance of one Gaussian from samples corrupted by noise, 
e.g. \cite{lai2016agnostic}, is run.
Once parameters of one Gaussian are estimated, all samples are ordered by Mahalanobis
distance \cite{mahalanobis1936generalised} to the Gaussian with the estimated mean and covariance.
The closest samples, which are most likely to be samples from this Gaussian with the estimated mean and covariance, 
are removed from future consideration and the algorithm proceeds with the estimation of the next Gaussian.

\begin{algorithm}[H]
\caption{Iterative Parameter Estimation in Noisy 2-GMM} 
\label{alg1}
\begin{algorithmic}[1]

\REQUIRE $X=\{ x_1, x_2,\ldots, x_m\}$, $w_1$.

\ENSURE $\hat{\mu}_1,\hat{\mu}_2$.

\STATE Let $\hat{\mu}_1=$AGNOSTICMEAN$(X)$. \\
       Let $\hat{\Sigma} =\text{AGNOSTICCOV}(X, 1-w_1)$.

\STATE For $i=1$ to $m$, $$y_i = (x_i -\hat{\mu}_1)^T \hat{\Sigma}^{-1} (x_i -\hat{\mu}_1).$$

\STATE Sort $$y_{(1)} \leq \ldots y_{([mw_1])} \leq \ldots y_{(m)}.$$ 

\STATE Let $X' = \{x_i: y_i >= y_{([mw_1])} \}$. 

\STATE $ \hat{\mu}_2 =\text{AGNOSTICMEAN}(X')$.

\end{algorithmic}
\end{algorithm}

See Algorithm~\ref{alg1} for the pseudo-code.
For completeness, we list the pseudo-code of algorithms AGNOSTICMEAN and AGNOSTICCOV of \cite{lai2016agnostic} in the appendix.

\section{The Sample Complexity}
\label{sec:complexity}

Our main result is the analysis of sample complexity of Algorithm~\ref{alg1}.

\begin{thm}[Sample Complexity of Iterative Parameter Estimation in Noisy 2-GMM]
\label{thm1}
For Problem~\ref{def:P1}, there exists a poly$(n,\frac{1}{\epsilon })$- time algorithm that takes as input $m$ independent 
samples $x_1,x_2,\ldots, x_m$ and $w_1$, then computes the means of two components of noisy 2-GMM model, s.t. \\
If $\Sigma = \sigma^2 I $, then
$$  \| \mu_2 -\hat{\mu}_2\|_2 = O({\frac{w_3}{w_2}} + \epsilon) \sigma \sqrt{ \log n}, $$
$$  \| \mu_1 -\hat{\mu}_1\|_2 = O(\frac{w_2}{w_1} + \epsilon)\sigma \sqrt{ \log n}. $$
Otherwise, for arbitrary $\Sigma$, 
$$  \| \mu_2 -\hat{\mu}_2\| = O(\sqrt{\frac{w_3}{w_2}} + \epsilon) \|\Sigma\|_2^{1/2} \sqrt{ \log n}, $$
$$  \| \mu_1 -\hat{\mu}_1\| = O(\sqrt{\frac{w_2}{w_1}} + \epsilon) \|\Sigma\|_2^{1/2} \sqrt{ \log n},  $$
provided that, 
\begin{enumerate}[i)]
\item  for the spherical case, $m= \Omega ( \frac{n \frac{1}{w_2}( \log n + \log \frac{1}{\epsilon}) \log n}{\epsilon^2} +  \frac{ \frac{1}{w_2} \log (\frac{1}{\epsilon} + \frac{w_2}{w_3}) }{\epsilon^2 + (\frac{w_3}{w_2})^2 } +  \frac{ \log (\frac{1}{\epsilon} + {\frac{w_2}{w_3} }) }{w_3} ) $.
\\ For the non-spherical case, $ m = \Omega ( \frac{n(n + \frac{1}{w_2})( \log n + \log \frac{1}{\epsilon}) \log n}{\epsilon^2} +   \frac{ \frac{1}{w_2} \log (\frac{1}{\epsilon} + \sqrt{\frac{w_2}{w_3} }) }{\epsilon^4 + (\frac{w_3}{w_2})^2 } +  \frac{ \log (\frac{1}{\epsilon} + \sqrt{\frac{w_2}{w_3} }) }{w_3} ).$ 
\item conditions in Lemma  \ref{lemma1nons}  are satisfied with $\eta = O(\epsilon + \frac{w_3}{w_2})$ in the spherical case $\Sigma = \sigma^2 I$, and $\eta = O(\epsilon^2 + \frac{w_3}{w_2})$ otherwise where $\eta$ is a parameter in Lemma \ref{lemma1nons}.
\end{enumerate}
\end{thm}

This result builds upon the analysis of  the separation of the two Gaussians.
Condition on recovery of the first mean comes from  \cite {lai2016agnostic}.
Our contribution lies in recovering the second component with mixing coefficient $w_2$ in a noisy 2-GMM, 
i.e., the degenerate component of the 2-GMM. 
In order to derive a robust estimator of $\mu_2$, we derive separation conditions on the two Gaussians in (ii) of Theorem \ref{thm1} 
with $\eta$ being the auxiliary parameter of accuracy of the estimator $\hat{\mu}_2$.
Under such conditions, one can filter and subsample according to Mahalanobis-distance criteria,
and take subsampled points as an input of the agnostic learning of $\mu_2$.

\subsection{Main Ideas of the Proof}

The complete proof is provided in the appendix. In this section, we show that there exists a polynomial-time 
algorithm that can estimate mean and covariance matrix of a Gaussian distribution from samples corrupted by 
malicious noise, i.e., an approximation algorithm for: 

\begin{defi}[Robust Parameter Estimation in Noisy 1-GMM]  \label{def:P2}
Given points in $\mathbb{R}^n$ that are each, with probability $1-\eta$ from an unknown 
distribution with mean $\mu$ and covariance $\Sigma$, and with probability $\eta$ completely arbitrary, estimate $\mu$ and $\Sigma$.
\end{defi}

Their most significant results are:

\begin{lemma}[Mean recovery in \cite{lai2016agnostic}]
For Problem~\ref{def:P2}, there exists a poly($n,1/\epsilon$)-time algorithm that takes as 
input $m =O ( \frac{n ( \log n  +\log \frac{1}{\epsilon} )\log n}{\epsilon^2} )$ independent 
samples $x_1,x_2,\ldots, x_m \sim N_{\eta}(\mu,\Sigma)$ and computes $\hat{\mu}$ such that
the error $\| \mu-\hat{\mu}\|_2$ is bounded as follows,
\begin{align*}
\begin{split}
O(\eta+ \epsilon) \sigma\sqrt{\log n}, & \quad \text{  if } \Sigma = \sigma^2 I. \\
O(\sqrt{\eta} + \epsilon) \|\Sigma\|_2^{1/2} \sqrt{\log n,} & \quad \text{ otherwise}.
\end{split}
\end{align*}
\label{thm0}
\end{lemma}

\begin{lemma}[Covariance recovery in \cite{lai2016agnostic}]
\label{thm01}
For Problem \ref{def:P2}, there exists a poly($n,1/\epsilon$)-time algorithm that takes as 
input $m = \Omega ( \frac{n^2 ( \log n  +\log \frac{1}{\epsilon} )\log n}{\epsilon^2} )$ independent samples
$x_1,x_2,\ldots, x_m \sim N_{\eta}(\mu,\Sigma)$ and computes $\hat{\Sigma}$ such that the error 
$$\| \Sigma-\hat{\Sigma}\|_F = O(\eta^{1/2} + C_1(\eta+\epsilon)^{3/4})\|\Sigma\|_2 \sqrt{\log n}$$
\end{lemma}

With the recovery of a Gaussian distribution from malicious noise, 
we can proceed with the case of imbalanced 2-GMM.
The crucial innovation is the condition on the separation between means of the 2-GMM
 based on the ratio of component weights and a careful analysis of the sample complexity.

\noindent \textbf{Notation.} 
Denote the set of the sampled points  $S=\{x_j\}_{j=1}^m$. $S={G_1} \cup {G_2} \cup N$, where ${G_k}$ denotes the samples from the $k^\text{th}$ component, and $N$ the set of  samples belonging to the malicious noise.
Denote by $$l(X) = (X- \hat{\mu}_1)'\hat{\Sigma}^{-1}(X-\hat{\mu}_1),$$
Mahalanobis distance. 
$l(X)_{(i )}$ is the $i^{\text{th}}$ smallest random variable among all samples $\{l(X_j)\}_{j=1}^m$.  
Let $l(X^{G_k})_{(i)}$ denote the $i^\text{th}$ smallest term among $m_k$ samples  from the $k^\text{th}$ component  $\{l(X_j): X_j \in G_k\}$.
Denote by
$$\lambda= (\mu_2- \mu_1)^T \Sigma^{-1}(\mu_2- \mu_1),$$
the `distance' of two Gaussian distribution in a 2-GMM. The larger $\lambda$ is, the better-separated the two Gaussian components are. 
In the spherical case $\Sigma= \sigma^2 I$, 
$$\lambda= \frac{\|\mu_2 - \mu_1\|_2^2}{\sigma^2},$$ which can be seen as a signal-to-noise parameter.
  Let $\tr(\Sigma)$ be the trace of the covariance matrix, and $\tr(\Sigma^2)$ the trace of its squared matrix.

\quad \smallskip

\noindent \textbf{Proof Sketch of Theorem \ref{thm1}.}
The identification of $\mu_1$ is trivial by Lemma \ref{thm0}. The challenge lies in  efficiently learning the distribution of the second component, 
which is solved by Step 4 and 5 in Algorithm \ref{alg1}. Formally, Step 4 (filter) can be translated as follows. Given $m$ i.i.d. samples $\x_1,\x_2,\ldots, \x_m$, 
we order the $l(X_i)$ by their magnitude,
$$l(X)_{(1)} \leq l(X)_{(2)} \leq  \ldots \leq l(X)_{(mw_1)} \leq \ldots \leq l(X)_{(m)},$$
and take as input $\{X_j : l(X_j) \geq  l(X)_{(mw_1)}\}$ in Step 5.

Succeeding in identifying samples from the second components leads to Theorem \ref{thm1}. In other words, if for any $\eta>0$, given sufficiently 
large sample size and two sufficiently well-separated Gaussian distributions in an imbalanced case,  
$$\Pr( w\in G_1 | l(w) \geq l(X)_{(mw_1)}) <\eta,$$
that is, among samples left over after Step 4 (Filter), those from the 1st component account for a sufficiently small proportion, then the input of Step 5 can be regarded as a Gaussian $N(\mu_2,\Sigma)$ with malicious noise of weight at most
$$\eta' = \eta + \frac{w_3}{w_2}.$$
Hence, applying the AGNOSTICMEAN algorithm in Step 5 can give a good estimator $\hat{\mu}_2$ of the second mean  given that $m (w_2+w_3)  = \Omega( \frac{n ( \log n  +\log \frac{1}{\epsilon} )\log n}{\epsilon^2} )$. With the success of step 4 and sufficiently large sample size we have,  the estimation error $\|\hat{\mu}_2 - \mu_2\|$  is bounded by 
\begin{align*}
\begin{split}
O(\eta + \frac{w_3}{w_2}+ \epsilon) \sigma\sqrt{\log n}, & \quad \text{  if } \Sigma = \sigma^2 I \\
O(\sqrt{\eta + \frac{w_3}{w_2}} + \epsilon) \|\Sigma\|_2^{1/2} \sqrt{\log n}, & \quad \text{ otherwise.}
\end{split}
\end{align*}
Let $\eta= O(\frac{w_3}{w_2}  + \epsilon)$ in the spherical case ($\Sigma = \sigma^2 I$), or $\eta=O( \epsilon^2 + \frac{w_3}{w_2})$ in the non-spherical case, we have the recovery result in Theorem \ref{thm1}.

Formally, Step 4 (filter) guarantees the following result.

\begin{lemma}[Non-Spherical Gaussian Mixture]
\label{lemma1nons}
For any $0<\eta<1$, given $m \geq \frac{ c w_2 \log \frac{w_2}{w_3}} {w_3 \eta^2}$ i.i.d. samples for some $c>0$,  we have
$$\Pr( w\in G_1 | l(w) \geq l(X)_{(mw_1)}) <\eta  ,$$ if the following conditions are satisfied
\begin{itemize}
\item $\lambda \geq 2 \sqrt{n\log\frac{w_1}{w_2}} + 2\sqrt{n-1} + 2\log\frac{w_1}{w_2}.$
\item Let $\delta=\eta \frac{w_2}{w_1} + \frac{w_3}{w_1}$, then
$$\lambda \geq  c_1 \sqrt{w_2\log n} ( \|\Sigma\|_2^{\frac{1}{2}} \|\Sigma^{-1}\|_2 +1) ( \tr(\Sigma)  +  \sqrt{\tr(\Sigma^2)\log\frac{1}{\delta}} + \|\Sigma\| \log\frac{1}{\delta} )  + c_2 \sqrt{n\log \frac{1}{\delta}}  + \log\frac{1}{\delta}.$$
If $n=o(\log \frac{1}{\delta})$,  it can be simplified as,
$$\lambda \geq \big(1+ c_1 \sqrt{w_2 \log n})( \|\Sigma\|_2^{\frac{1}{2}} \|\Sigma^{-1}\|_2 +1) \|\Sigma\|_2  \log\frac{1}{\delta} .$$
\item The smallest singular value of $\Sigma$ is bounded away from $0$, i.e.,
$\sigma_{\min}=\| \Sigma^{-1} \|_2^{-2} > O(w_2)\|\Sigma\|_2 \log n.$\\
 $\sigma_{\min} \geq O\big(\sqrt{w_2}  \|\Sigma\|_2 \sqrt{\log n} (\frac{  2\|\Sigma\|_2^2 \sqrt{\log \frac{w_2}{w_3}} }{ \sqrt{\tr(\Sigma^2)} } +1) \big)$ 
\end{itemize}
\end{lemma}

In the spherical case, the filtering guarantee can be simplified.

\begin{lemma}[Spherical Gaussian Mixture]
\label{lemma1s}
For any $0<\eta<1$, with $m \geq \frac{ c w_2 \log \frac{w_2}{w_3}} {w_3 \eta^2}$ i.i.d. samples for some $c>0$, 
$$\Pr( w\in G_1 | l(w) \geq l(X)_{(m_1)}) <\eta  ,$$ if the following conditions are satisfied
\begin{itemize}
\item $\lambda \geq 2 \sqrt{n\log\frac{w_1}{w_2}} + 2\sqrt{n-1} + 2\log\frac{w_1}{w_2}.$
\item Let $\delta=\eta \frac{w_2}{w_1} + \frac{w_3}{w_1}$, then 
  $$\lambda \geq   c_1 \log{\frac{1}{\delta}} + c_2\sqrt{n\log{\frac{1}{\delta} } } + c_3 \sqrt{\log\frac{1}{\eta}( \log\frac{1}{\delta} +n )} ,$$
for some $c_1, c_2, c_3 >0$.\\
If $n=o(\log \frac{1}{\delta})$, then it is equivalent to,
$$\lambda \geq c'_1 \log\frac{1}{\delta}, $$
for some $c'_1$.
\end{itemize}
\end{lemma}

\noindent \textbf{Remark 1.} The separation of $\lambda$ depends  on three parameters of the 2-GMM models. 
\begin{enumerate}
\item The ratio of two components $\frac{w_2}{w_1}$.
\item The accuracy of estimation $\delta= \eta\frac{w_2}{w_1}  + \frac{w_3}{w_1}$.
\item The dimension of the problem $n$.
\end{enumerate}
For $\frac{w_2}{w_1}$, it is obvious that the more skewed (imbalanced) the 2-GMM model is, the more difficult it is to learn the means of the two components efficiently. Secondly, notice that $\delta=\frac{w_2}{w_1}(\eta+ \frac{w_3}{w_2})$, thus the order of $\delta$ is between the order of the two ratios. $\frac{w_2}{w_1}$ is the generic upper bound on the accuracy of estimation, while $\eta$ is an auxiliary parameter for the accuracy of estimation one would like to achieve in estimating the smaller component.
Therefore, the more accurate the estimation is for the second mean, i.e., the smaller $\eta$ is, the more strict separation conditions on $\lambda$ are.
On the other hand, the accuracy term $\delta$ is bounded below  by $\frac{w_3}{w_1}$ indicating that the underlying bound posed by the malicious noise. Thirdly, separation conditions between the two components depend on the dimension of the problem in an imbalanced case. In particular, it requires that in each one-dimensional direction, the mean of skewed distributed component $\mu_2$ is roughly $O(n^{-1/4})$ away from $\mu_1$ in that direction. This is a stronger separation condition compared to the balanced case (\cite{balakrishnan2017statistical}).

\smallskip

\noindent \textbf{Remark 2.}  In the non-spherical case, the third condition in Lemma \ref{lemma1nons} can be translated into $\frac{\sigma_{\min}}{\sigma_{\max}}$ being bounded away from $0$.  This contributes to a good approximation of $\Sigma^{-1}$ using $\hat{\Sigma}$ which is close to $\Sigma$ in Frobenius norm.

The following is a proof sketch of the crucial filtering result in the spherical case. The proof for the non-spherical case follows analogously. 

\bigskip

\noindent \textbf{Proof Sketch of Lemma \ref{lemma1s} }
Denote by $m_1 = w_1m$. By Bayes rule,  
\begin{align*}
\begin{split}
& \Pr( w \in G_1 | l(w) \geq l(X)_{(m_1 )})  \\
 = & \frac{\Pr( l(w) \geq l(X)_{(m_1)} | w \in G_1 ) \Pr( w \in G_1) }{ \sum_{J =G_1, G_2, N}  \Pr(  l(w) \geq l(X)_{(m_1)} | w \in J) \Pr( w \in J ) } \\
 \leq & \frac{\Pr( l(w) \geq l(X)_{(m_1)} | w \in G_1 ) \Pr( w \in G_1) }{  \sum_{J =G_1, G_2}  \Pr(  l(w) \geq l(X)_{(m_1)} | w \in J) \Pr( w \in J )   }  \\
 = & \frac{\Pr( l(w) \geq l(X)_{(m_1)} | w \in G_1 ) }{\Pr( l(w) \geq l(X)_{(m_1)} | w \in G_1 )  + \frac{w_2}{w_1} \Pr( l(w) \geq l(X)_{(m_1)} | w \in G_2 ) }.
  \end{split}
 \end{align*}
 Therefore, it suffices to show that,
\begin{align}
\begin{split}
\label{eq01}
& \Pr( l(w) \geq l(X)_{(m_1)} | w \in G_1 ) \\
\leq & \eta  \frac{w_2}{w_1} \Pr( l(w) \geq l(X)_{(m_1)} | w \in G_2 ). 
\end{split}
\end{align}
To show \eqref{eq01}, we show an upper bound of LHS and lower bound of RHS and find conditions so that the upper bound  of LHS is smaller than the lower bound of RHS. 
The proof can be decomposed into 2 steps.

\noindent \textbf{(i)} An upper bound of $\Pr( l(w) \geq l(X)_{(m_1)} | w \in G_1 )$
\begin{align*}
\begin{split}
 & \Pr( l(w) \geq l(X)_{(m_1)} | w \in G_1 )  \\
 \leq & \Pr( l(w) \geq l(X^{G_1\cup G_2})_{(m_1-m_3)} | w \in G_1  ) \\
  = &  \Pr(l(w) \geq l(X^{G_1\cup G_2} )_{(m_1-m_3)}, l(w) \leq l(X^{G_1})_{(m_1(1-\beta) -m_3 )} | w \in G_1  ) \\ 
   + & \Pr(l(w) \geq l(X^{G_1\cup G_2} )_{(m_1-m_3)}, l(w) \geq l(X^{G_1})_{(m_1(1-\beta) -m_3) } | w \in G_1  ) \\
  \leq & \sum_{n_1+n_2 = m_1 - m_3}\sum_{j=n_1}^{m_1(1-\beta)-m_3}  \Pr(  l(X^{G_1})_{(j) }=l(w) \geq l(X^{G_2})_{(n_2)} |w\in G_1) + \beta +\frac{m_3}{m_1}.
 \end{split}
 \end{align*}
 Notice that $\Pr(  l(X^{G_1})_{(j) } \geq l(X^{G_2})_{(n_2)}) $ is bounded by
 \begin{align*}
 \begin{split}  
 & \Pr \big(  l(X^{G_1})_{(m_1(1-\beta) -m_3) } \geq l(X^{G_2})_{(m_1\beta ) } \big)  \\
 \leq & \Pr\big(l(X^{G_1})_{(m_1(1-\beta) -m_3) }\geq t)  + \Pr( l(X^{G_2})_{(m_1\beta) } \leq t)\\
  \leq &\Pr \left(\frac{S_{m_1}(\underline{p})}{m_1} \leq 1-\beta-\frac{m_3}{m_1} \right) + \Pr \left(\frac{S_{m_2}(\bar{p})}{m_2} \geq \frac{m_1 \beta}{m_2} \right),
 \end{split}
 \end{align*}
for some $t,\bar{p},\underline{p}$ such that $\underline{p} \leq \Pr(l(X_i^{G_1}) \leq t)$ and $\bar{p} \geq \Pr(l(X_i^{G_2}) \leq t) $ and $S_n(p)$ denotes the sum of $n$ i.i.d. Bernoulli trials with parameter $p$.  
The reason of introducing $t$ is that $l(X_i^{G_2})$ is approximately (under $\hat{\mu}_1, \hat{\Sigma}_1$)  non-central chi-squared distributed.\footnote{A non-central chi-squared distribution is denoted as $\chi^2_n(\lambda)$ with $\lambda$ being the non-centrality parameter and $n$ the degree of freedom. Let $(Y_1, Y_2, \ldots, Y_n)$ be $n$ independent normally distributed random variables with mean $a_i$ and unit variance. Then $\sum_{j=1}^{n} Y_j^2 $ follows a $\chi^2_n(\lambda)$ distribution and $\lambda =\sum_{j} a_j^2.$  }
With a well-chosen $t$ and the condition that the two Gaussians are well-separated,  we can approximate the comparison between a central and a non-central chi-squared distribution by finding a cutoff $t$ that well separates the two ellipsoids, instead of working out a joint distribution of the two.

 For the first term and $\underline{p} \leq \Pr(l(X_i^{G_1}) \leq t)$, we apply a sharp bound on the right tail of central chi-squared distribution \cite{laurent2000adaptive},
 $$ \Pr( \chi^2_n -n \geq 2\sqrt{xn} + 2x ) \leq \exp(-x).$$
Subsequently, we set $\underline{p}$ accordingly.

For the second term and $\bar{p} \geq \Pr(l(X_i^{G_2}) \leq t) $,  we apply similarly an upper bound on the left tail of non-central chi-squared distribution $\chi_n^2(\lambda)$  given by \cite{birge2001alternative},
$$\Pr(\chi_n^2(\lambda)  \leq (n+\lambda) - 2\sqrt{(n+2\lambda)x} ) \leq \exp(-x).$$
which determines our $\bar{p}$ accordingly.

In the proof, we show that under the well-separated condition on $\lambda$, one could find satisfying threshold $t$, such that, 
$$\Pr \left(\frac{S_{m_1}(\underline{p})}{m_1} \leq 1-\beta-\frac{m_3}{m_1} \right) = O(\exp( - m_1 \beta)),$$
$$\Pr \left(\frac{S_{m_2}(\bar{p})}{m_2} \geq \frac{m_1 \beta}{m_2} \right)= O(\exp(- m_1 \beta)). $$
Hence, 
$$\Pr(  l(X^{G_1})_{(m_1(1-\beta) -m_3) } \geq l(X^{G_2})_{(m_1\beta )}) = O(\exp(- m_1 \beta)).$$
 Therefore, under the choice of $\beta = \frac{m_2}{m_1}\eta ,$ 
\begin{align*}
\begin{split}
\Pr( l(w) \geq l(X)_{(m_1)} | w \in G_1 ) = O(\frac{w_2}{w_1} \eta + \frac{w_3}{w_1} ).
\end{split}
\end{align*}

\noindent \textbf{(ii)} A  lower bound on $\Pr(  l(w) \geq l(X)_{(m_1)} | w \in G_2)$.

 To quantify the right tail of a non-central chi-squared distribution, we  apply an alternative characterization  of the $\chi^2_n(\lambda)$ as follows.
\begin{lemma}[\cite{cacoullos1984quadratic}]
\label{lemma: chisquare1}
$$\chi_n^2(\lambda) \myeq \chi_{n-1}^2(0) +  \chi_1^2(\lambda).$$
\end{lemma}
In other words, a non-central chi-squared distribution with degree of freedom $n$ and non-centrality parameter $\lambda$ is, in distribution, equivalent to the sum of a random variable drawn from non-central chi-squared distribution with degree of freedom $1$ and non-centrality parameter $\lambda$ and another random variable drawn from central chi-squared distribution with degree of freedom $n-1$. Moreover, the  two random variables are independent.  

Then a lower bound on $\Pr(  l(w) \geq l(X)_{(m_1)} | w \in G_2)$ can be obtained with a proper choice of $\tilde{m}_2=O( m_2)$ such that, 
\begin{align*}
\begin{split}
& \Pr( l(w) \geq l(X)_{(m_1)} | w\in G_2) \\
\geq &  \sum_{i + j= m_1} \Pr(l(w) \geq l(X^{G_1})_{(i)}, l(w) \geq l(X^{G_2})_{(j)} | w\in G_2) \\
 \geq & \Pr(l(w) > t , l(w) \geq l(X^{G_2})_{(\tilde{m}_2 )} | w\in G_2 )  \Pr(  l(X^{G_1})_{(m_1 - \tilde{m}_2 )} \leq t ) \\ 
 = & \frac{m_2 - \tilde{m}_2}{m_2} \Pr(l(X^{G_2})_{(\tilde{m}_2 )} > t )  \Pr(  l(X^{G_1})_{(m_1 - \tilde{m}_2 )} \leq t ).
\end{split}
\end{align*}
Similarly we would like to find a ball of `radius' $t$ that covers at least $1-\frac{\tilde{m}_2}{m_1}$ fraction of sampled points from the first component, while  overlapping at most a fraction $\frac{\tilde{m}_2}{m_2}$ of the samples from the second components.   Using a tail bound for $\chi^2_n$ given by  \cite{laurent2000adaptive},
one could prove that $$\Pr(  l(X^{G_1})_{(m_1 - \tilde{m}_2 )} \leq t ) \geq  1- c \exp ( -\tilde{m}_2 ) .    $$
On the other hand, under the well-separated condition, and one could show that  with the choice of $\tilde{m}_2$,
$$ \Pr(l(X^{G_2})_{(m_2-\tilde{m}_2 : m_2)} > t ) \geq 1/2.$$
Therefore, for the choice of $\tilde{m}_2$ to be specified in the appendix, $\Pr( l(w) \geq l(X)_{(m_1)} | w\in G_2)$ is bounded below by some constant when $m_2$ is sufficiently large.

One could show that under the separation condition, the upper bound of LHS in \eqref{eq01} is smaller than the lower bound of its RHS, thus completing the proof of  Lemma \ref{lemma1s}. The difference between the proofs for the spherical and the non-spherical case is resolved by \cite{hsu2012tail} with the following inequalities.  If $y\sim N(\mu,\Sigma)$, then
$$\Pr(\|y-\mu\|^2 \geq \tr(\Sigma) + 2 \sqrt{\tr(\Sigma^2)x} + 2\|\Sigma\|x )\leq \exp(-x).$$
Moreover, if $y\sim N(\mu_2,\Sigma)$,
\begin{align*}
\begin{split}
\Pr(\|y-\mu_1\|^2\geq \tr(\Sigma) + 2 \sqrt{\tr(\Sigma^2) x} + 2 \|\Sigma\| x +  \|\mu_2 - \mu_1\|^2 (1+\frac{2\|\Sigma \| x}{\sqrt{\tr(\Sigma^2)x}}) \leq \exp(-x).
\end{split}
\end{align*}
Then, applying the same idea of find a ball of `radius' $t$ for some proper choice of $t$ that separates the two components, we can achieve the desired inequalities.

\section{Sensitivity}
\label{sec:sensitivity}

Compared with the widely adopted \cite{titterington1985statistical} vanilla expectation-maximisation (EM) algorithm, 
our algorithm for recovering the mean of 2-GMMs does not requires an initialisation. 
Moreover, the algorithm significantly outperforms the EM algorithm given the same number of sampled points,
as demonstrated in the simulations of the next section. 
Nevertheless, it is worth noticing that our algorithm does require the mixing coefficient $w_1$ on the input. 
In the spherical case, we show that the algorithm is robust to a perturbation in the  input $w_1$.
That is, if instead of $w_1$, the input takes an imprecise estimator $w'_1$, the output of the estimated means $\hat{\mu}_1, \hat{\mu}_2$ are slightly perturbed.
\begin{prop}
\label{prop_sensitive}
Consider the spherical case and assume that conditions of Theorem \ref{thm1} are satisfied. Then, there exist $\alpha_k, k=1,2, \ldots$ such that 
$\forall 1>w'_1>0$ we have $\|\hat{\mu}_k (w'_1) - \hat{\mu}_k (w_1)\|_2$ bounded by:
\begin{equation*}
\|\hat{\mu}_k (w'_1) - \hat{\mu}_k (w_1)\|_2 \leq \alpha_k (\frac{|w'_1 - w_1| + w_3}{w_2}  + \frac{w_3}{1-w'_1} +  {\epsilon }) \sigma \sqrt{\log n }.
\end{equation*}
\end{prop}

See Figures~\ref{sen1} and \ref{sen2} in the next section for a computational illustration.

\section{Computational Illustrations}
\label{sec:computational}

We performed 50 simulations on $10,000$ samples from a noisy 2-GMM each, with mixing coefficients $w_1= 0.8,w_2=0.16, w_3=0.04$ throughout. 
The noise is Cauchy distributed and dimensions $n= 10,12,14,16,18,$ $20,40,60,80$. 
Table \ref{table1} compares the estimation error on $\mu_1$ and $\mu_2$ achieved by Algorithm \ref{alg1} to the estimation error achieved by the vanilla EM algorithm,
implemented as fitgmdist in Mathworks Matlab 2016a.
Following the conventional rules, we measure the estimation error as
$$   \min_{\pi \in P_2} \|\hat{\mu}_{\pi(1)} - \mu_1\|_2  + \|\hat{\mu}_{\pi(2)} - \mu_1\|_2.$$
The true sampling error in the graph denotes the labeled mean $\|mean(X^{G_1}) - \mu_1\|$, $\|mean(X^{G_2}) - \mu_2\|$.
We start with the 2-GMM with malicious noise in $n= 10,12,14,16,18,$ $20$ dimensional  cases.

\begin{table*}[t]
\centering
\small
\begin{tabular}{cl|llllll}
Estimation error                                                   & \multicolumn{1}{c|}{} & \multicolumn{1}{c}{$n=10$} & \multicolumn{1}{c}{$n=12$} & \multicolumn{1}{c}{$n=14$} & \multicolumn{1}{c}{$n=16$} & \multicolumn{1}{c}{$n=18$} & \multicolumn{1}{c}{$n=20$}  \\ \hline
\multirow{2}{*}{$\|\hat{\mu}_1- \mu_1\|$}                          & Alg. \ref{alg1}               & 0.54(0.13)             & 0.51(0.10)             & 0.53(0.13)             & 0.56(0.14)             & 0.55(0.10)             & 0.55(0.13)        \\
                                                                   & EM                     & 1.66(2.56)             & 2.04(3.32)             & 1.69(1.61)             & 1.69(1.21)             & 2.08(2.41)             & 1.72(0.52)    \\ \hline
\multirow{2}{*}{$\|\hat{\mu}_2 - \mu_2\|$}                         & Alg. \ref{alg1}               & 1.18(0.19)             & 1.17(0.13)             & 1.16(0.19)             & 1.22(0.17)             & 1.25(0.16)             & 1.26(0.16)        \\
                                                                   & EM                     & 78.71(203.58)          & 21.07(15.34)           & 50.21(62.11)           & 78.04(276.22)          & 39.97(60.06)           & 60.27(97.50)     \\ \hline
\multirow{2}{*}{$\|\hat{\mu}_1 - \mu_1\|+\|\hat{\mu}_2 - \mu_2\|$} & Alg. \ref{alg1}               & 1.72(0.27)             & 1.68(0.16)             & 1.70(0.25)             & 1.78(0.26)             & 1.80(0.19)             & 1.81(0.22)    \\
                                                                   & EM                     & 80.37(203.45)          & 23.11(14.97)           & 51.90(62.16)           & 79.73(276.13)          & 42.05(59.86)           & 61.99(97.45)         
\end{tabular}
\caption{Estimation errors on $\mu_1$ and $\mu_2$ using Algorithm~\ref{alg1} and Vanilla EM}
\label{table1}
\end{table*}

\begin{table*}[t]
\centering
\small
\begin{tabular}{cl|lllll}
Estimation error                                                   & \multicolumn{1}{c|}{} & \multicolumn{1}{c}{$n=10$} & \multicolumn{1}{c}{$n=20$} & \multicolumn{1}{c}{$n=40$} & \multicolumn{1}{c}{$n=60$}  & \multicolumn{1}{c}{$n=80$} \\ \hline
\multirow{2}{*}{$\|\hat{\mu}_1- \mu_1\|$}                          & Alg. \ref{alg1}               & 0.54(0.13)                      & 0.55(0.13)        & 0.77(0.26)     & 0.74(0.32)    & 0.77(0.33)   \\
                                                                   & EM                     & 1.66(2.56)                & 1.72(0.52)    &  2.99(2.46)  & 11.04(43.41)   &    7.95(11.42)    \\ \hline
\multirow{2}{*}{$\|\hat{\mu}_2 - \mu_2\|$}                         & Alg. \ref{alg1}               & 1.18(0.19)                    & 1.26(0.16)   & 1.59(0.20) & 3.54(0.81) & 1.88(0.19)          \\
                                                                   & EM                     & 78.71(203.58)               & 60.27(97.50)     & 156.80(491.43) & 221.60(040.94) & 98.24(161.98)      \\ \hline
\multirow{2}{*}{$\|\hat{\mu}_1 - \mu_1\|+\|\hat{\mu}_2 - \mu_2\|$} & Alg. \ref{alg1}               & 1.72(0.27)                      & 1.81(0.22)      & 2.37(0.35)   & 2.47(0.31)    & 2.64(0.44) \\
                                                                   & EM                     & 80.37(203.45)           & 61.99(97.45)          & 159.78(491.24)  & 232.64(1043.11) & 106.20(161.84)
\end{tabular}
\caption{Estimation errors on $\mu_1$ and $\mu_2$ using Algorithm~\ref{alg1} and Vanilla EM on Higher Dimensions}
\label{table2}
\end{table*}
\begin{table*}[t]
\centering
\begin{tabular}{c|llllll}
             & \multicolumn{1}{c}{$n=10$} & \multicolumn{1}{c}{$n=12$} & \multicolumn{1}{c}{$n=14$} & \multicolumn{1}{c}{$n=16$} & \multicolumn{1}{c}{$n=18$} & \multicolumn{1}{c}{$n=20$} \\ \hline
$\alpha = 0.50$  & 3.05(0.40)             & 3.09(0.33)             & 3.22(0.44)             & 3.18(0.37)             & 3.28(0.40)             & 3.39(0.44)             \\
$\alpha = 0.67$ & 2.65(0.40)             & 2.72(0.28)             & 2.82(0.41)             & 2.89(0.48)             & 2.90(0.34)             & 2.97(0.44)             \\
$\alpha = 0.83$ & 2.24(0.36)             & 2.25(0.27)             & 2.31(0.32)             & 2.39(0.37)             & 2.41(0.29)             & 2.53(0.33)             \\
$\alpha = 1.00$    & 1.72(0.27)             & 1.68(0.16)             & 1.70(0.25)             & 1.78(0.26)             & 1.80(0.18)             & 1.81(0.22)             \\
$\alpha = 1.17$ & 1.83(0.29)             & 1.78(0.20)             & 1.84(0.26)             & 1.95(0.29)             & 1.98(0.18)             & 2.05(0.24)             \\
$\alpha = 1.33$ & 2.05(0.32)             & 2.05(0.23)             & 2.16(0.30)             & 2.30(0.30)             & 2.31(0.24)             & 2.38(0.26)             \\
$\alpha = 1.50$  & 2.43(0.36)             & 2.53(0.29)             & 2.69(0.30)             & 2.82(0.35)             & 2.87(0.30)             & 3.00(0.32)         
\end{tabular}
\caption{Estimation error of $\|\hat{\mu}_1 - \mu_1\|+\|\hat{\mu}_2 - \mu_2\|$ as $\alpha = \frac{1-w'_1}{1-w_1}$ varies for higher dimensions }   
\label{my-label}
\end{table*}

\begin{figure}[th]
\centering
        \includegraphics[width=0.6\textwidth]{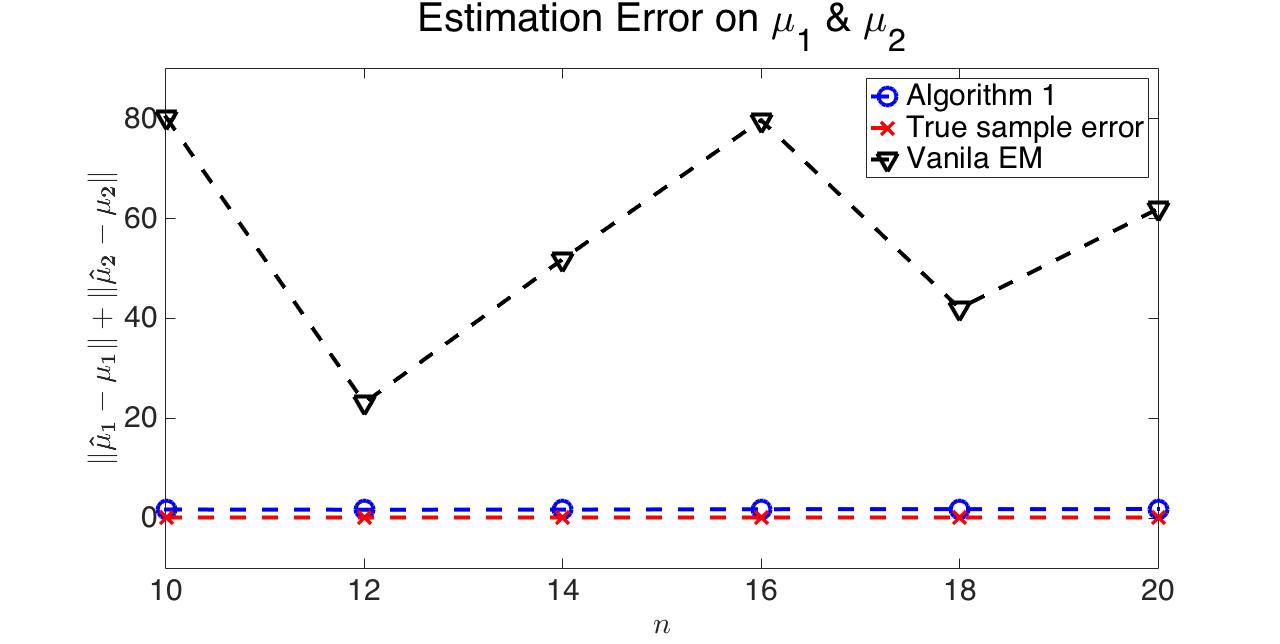} 
                \caption{A comparison of vanilla EM with Algorithm \ref{alg1}  in terms of estimation error $\|\hat{\mu}_1 -\mu_1\|_2+\|\hat{\mu}_2 -\mu_2\|_2$}
                \label{1g2g}
\end{figure}

\begin{figure}[th]
        \centering
        \includegraphics[width=0.6\textwidth]{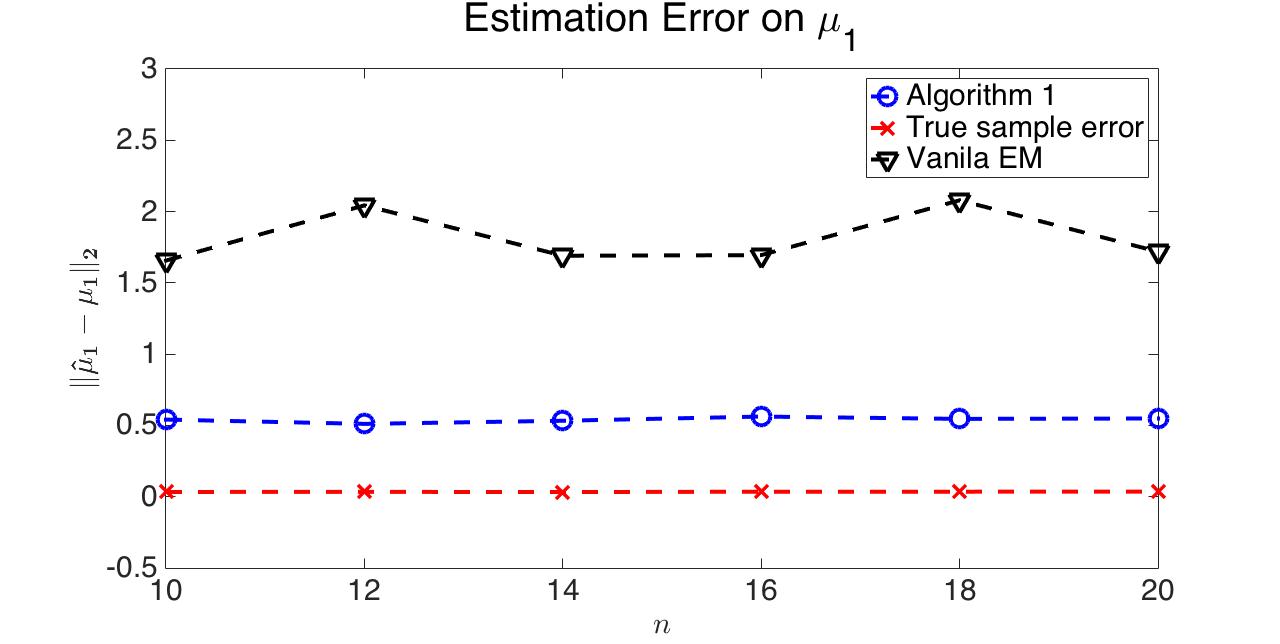} 
                        \caption{A comparison of vanilla EM with Algorithm \ref{alg1}  in terms of estimation error $\|\hat{\mu}_1-\mu_1\|_2$}
                        \label{1g}
\end{figure}

\begin{figure}[th]
        \centering
        \includegraphics[width=0.6\textwidth]{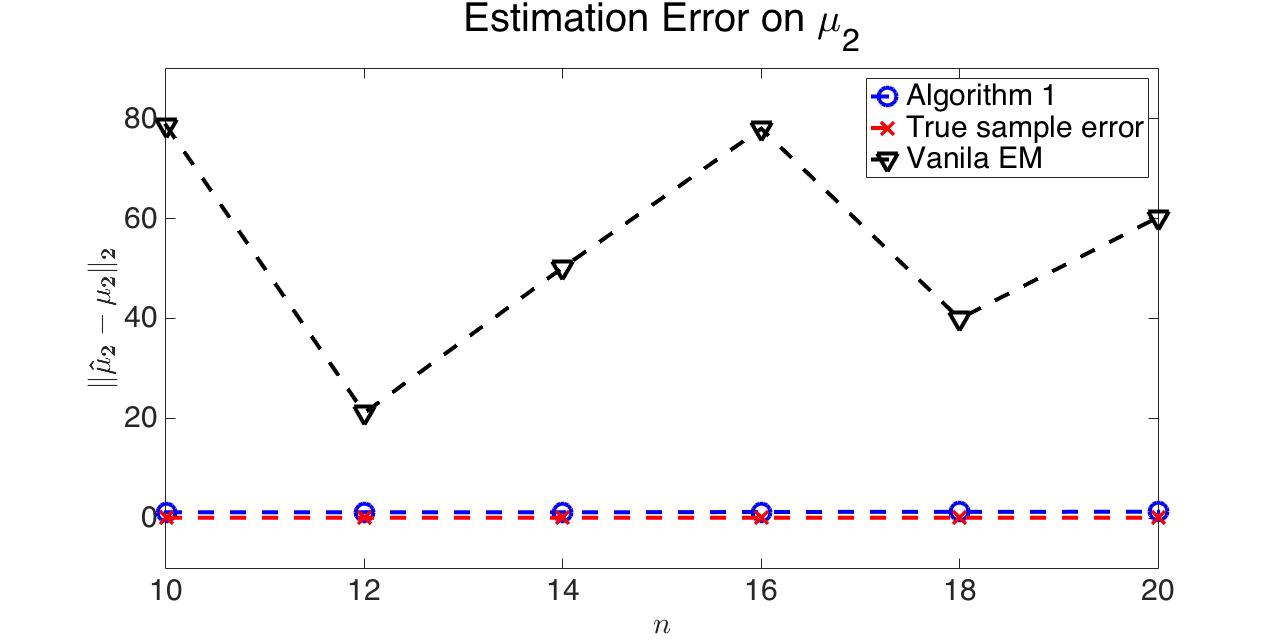} 
                        \caption{A comparison of vanilla EM with Algorithm \ref{alg1} in terms of estimation error $\|\hat{\mu}_2 -\mu_2\|_2$}
    \label{2g}
\end{figure}

Figures \ref{1g},\ref{2g},\ref{1g2g} show that the estimation errors by vanilla EM are much larger than that by Algorithm \ref{alg1}, 
especially in the second component. Based on Table \ref{table1}, the variance of the predictions by Algorithm \ref{alg1} is also  smaller than that of EM.

\begin{figure}[th]
\centering
        \includegraphics[width=0.6\textwidth]{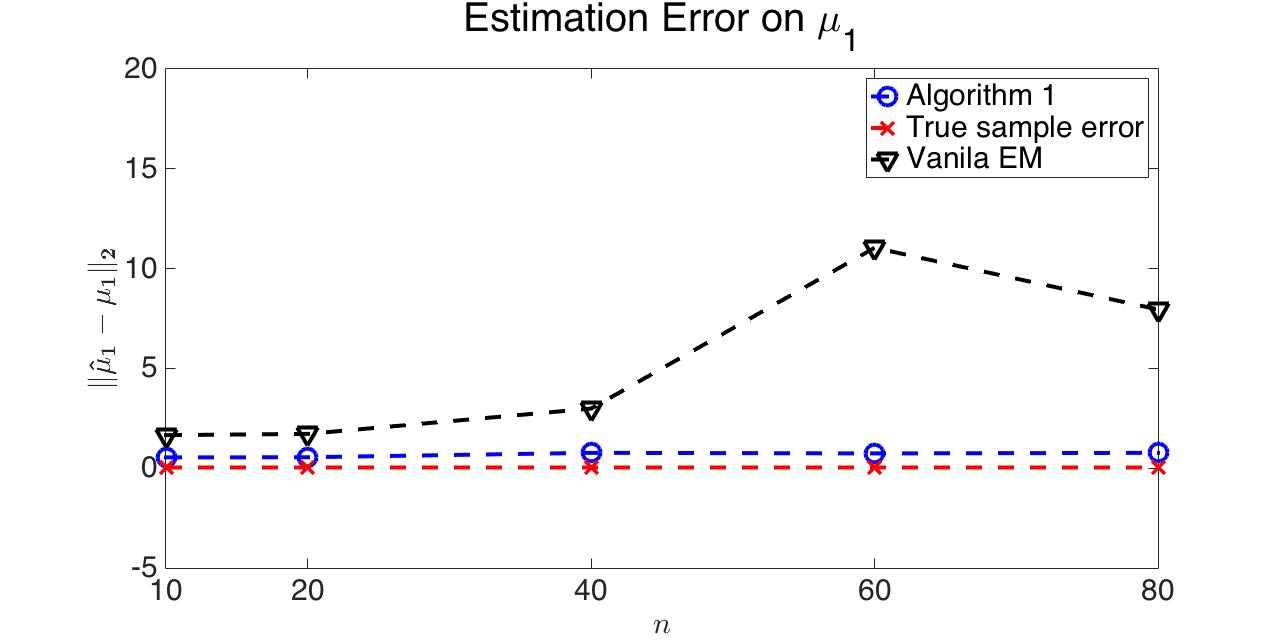} 
                \caption{A comparison of vanilla EM with Algorithm \ref{alg1} in terms of estimation error  $\|\hat{\mu}_1 -\mu_1\|_2 $ }
\label{high1}
\end{figure}

\begin{figure}[th]
        \centering
        \includegraphics[width=0.6\textwidth]{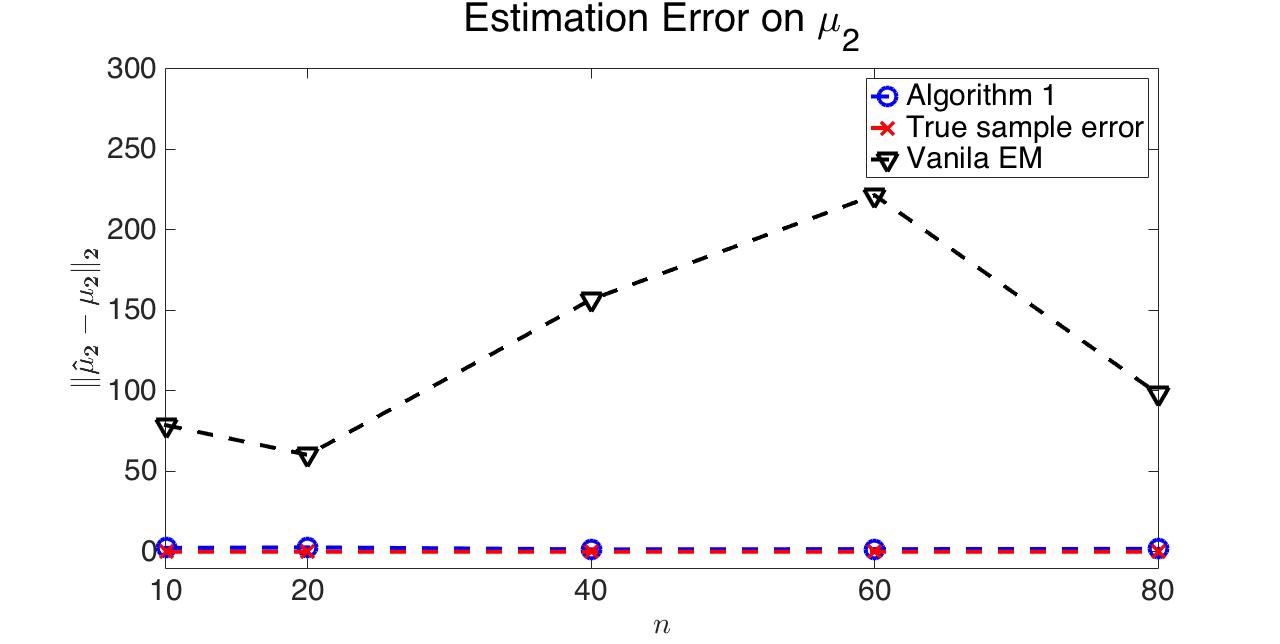} 
                        \caption{A comparison of vanilla EM with Algorithm \ref{alg1}  in terms of estimation error  $\|\hat{\mu}_2 -\mu_2\|_2 $}
\label{high2}
\end{figure}

From Figures \ref{high1},\ref{high2}, one could notice that as dimension goes up, vanilla EM algorithm witnesses much larger error and variance, compared with the more robust Algorithm \ref{alg1}, which does not requires any initialization.

\begin{figure}[th]
\centering
        \includegraphics[width=0.6\textwidth]{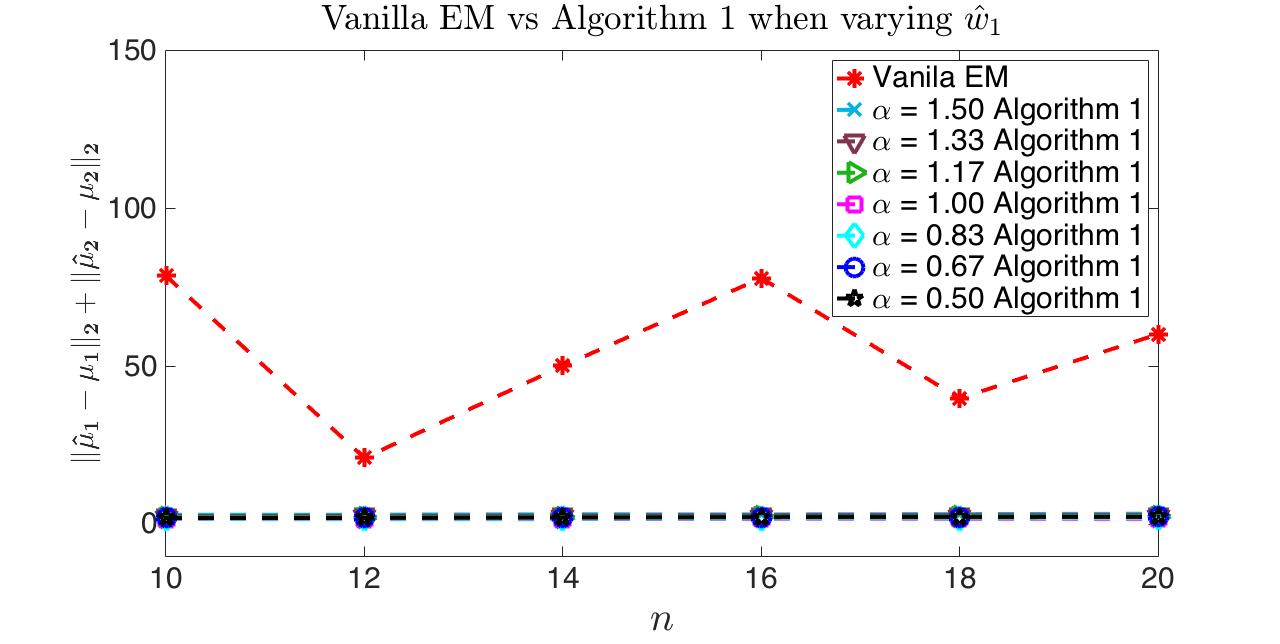} 
                \caption{A comparison of vanilla EM with Algorithm \ref{alg1} in terms of estimation error  $\|\hat{\mu}_1 -\mu_1\|_2 + \|\hat{\mu}_2-\mu_2\|_2$ }
\label{sen1}
\end{figure}

\begin{figure}[th]
        \centering
        \includegraphics[width=0.6\textwidth]{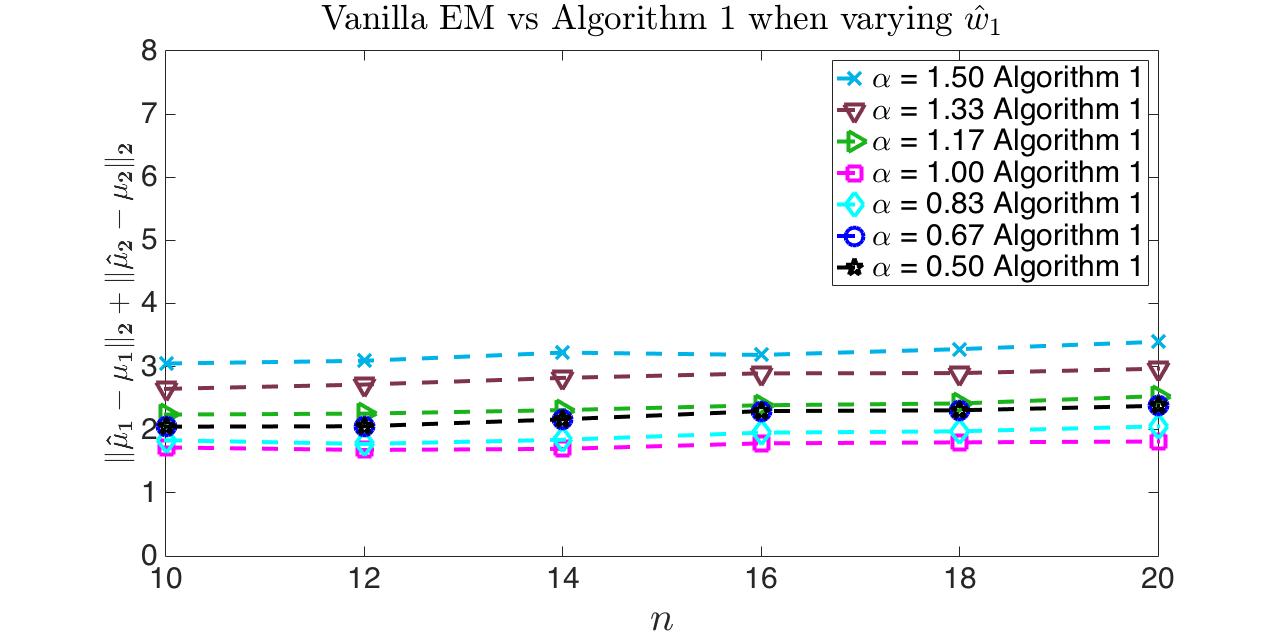} 
                        \caption{A comparison of vanilla EM with Algorithm \ref{alg1}  in terms of estimation error  $\|\hat{\mu}_1 -\mu_1\|_2 + \|\hat{\mu}_2-\mu_2\|_2$}
\label{sen2}
\end{figure}

In Figures \ref{2g} and \ref{1g2g}, while the average $\|\hat{\mu}_2 - \mu_2\|$ 
obtained by Algorithm~\ref{alg1} is below $5$, the averaged estimation error for the second component is above 20 for vanilla EM. 
From Table \ref{table1} we can further notice that
the outputs of the vanilla EM algorithm also suffer from large variation both compared with Algorithm \ref{alg1} as well as across different dimensions, which is due to its sensitivity to the initialization points.

Next, we conduct a more detailed experiment on the sensitivity to the input of $w_1$ in Algorithm~\ref{alg1}. 
Figures \ref{sen1} and \ref{sen2} show the estimation error on $\mu_1, \mu_2$ when varying 
$$\alpha = \frac{1-w'_1}{1-w_1}.$$
Clearly, we overestimate the mixing coefficient of the larger component, whenever $\alpha<1$.
For $\alpha>1$, the proportion of the smaller component and malicious noise, i.e. $w_1$, is overestimated. 
Figures~\ref{sen1} and \ref{sen2} outline the results, suggesting the performance of Algorithm~\ref{alg1} 
when the input $w'_1$ is distorted.

\section{Discussion}

We have presented a meta-algorithm for Robust Parameter Estimation in Noisy 2-GMM, which 
works best when the differences in mixing coefficients are large, i.e., $w_1 \gg w_2 \gg w_3$,
and the coefficients are known exactly.
There,
it outperforms the widely-used expectation-maximisation (EM) algorithm considerably, 
as documented in Section~\ref{sec:computational}.
The algorithm does not require any initialisation and it is rather stable with respect to the error in the estimate of the mixing coefficients,
as detailed in Section~\ref{sec:sensitivity}.
Our main result is an analysis of the sample complexity of the algorithm,
utilising spectral methods of \cite{lai2016agnostic}.

To continue in this direction, one could parametrise our analysis by the performance of the algorithm for Problem~\ref{def:P2},
and analyse whether EM algorithms for Problem~\ref{def:P2} would improve the performance.
Using information-theoretic arguments \cite{Yang1999}, one could also study the tightness of the bound on the sample complexity.  

\FloatBarrier
\bibliographystyle{ieeetr}
\bibliography{myreference}

\appendix

\onecolumn
\section{Algorithms }

\begin{algorithm}[H]
\caption{OUTLIERDAMPING(S)}
\begin{algorithmic}[1]

\REQUIRE $S=\{ \x_1, \x_2,\ldots, \x_m\} \subset \mathbb{R}^n $.

\ENSURE $S \subset \mathbb{R}^n, \mathbf{\omega}=(w_1,\ldots,w_m)$.

\STATE \textbf{If}  $n=1$:  \\
\textbf{Return} $(S,-1)$

\STATE Let $\mathbf{a}$ be the coordinate-wise media of $X$. Let $s^2 = C \text{Tr}(\Sigma)$. Estimate $\tr(\Sigma)$ by estimating $1$d variance along $n$ orthogonal directions and adding up.\footnotemark{}

\STATE Set $w_i = \exp(-\frac{\|x_i - \mathbf{a}\|_2^2}{s^2} )$ for every $\x_i \in X$.

\STATE \textbf{Return} $(S,\mathbf{\omega})$
\end{algorithmic}
\end{algorithm}

\footnotetext{ {}{Suppose $D = N(\mu,\sigma^2)$. Denote by  $D_{\eta}$  the 1-GMM with malicious noise of fraction $\eta$. There are several ways to estimate $\sigma$. One is, let $x_{med} = \text{MEDIAN}\{x_i\}$. Let $\Phi(x)$ be the c.d.f. of $N(0,1)$. Note that $c_1 = \Phi(1) \approx 0.851$. Let $C_{\sigma}$ be the $c_1 \times 100\%$th quantile of $S$. Let the estimator $\hat{\sigma} = C- \hat{\mu}$ where $\hat{\mu}$ is the estimated mean. One can prove that $|\hat{\sigma}^2 - \sigma^2| = O(\eta\sigma^2)$. } }

\begin{algorithm}[H]
\caption{OUTLIERTRUNCATION($S,\eta,\epsilon$)}
\begin{algorithmic}[1]

\REQUIRE $S=\{ \x_1, \x_2,\ldots, \x_m\} \subset \mathbb{R}^n $, $\eta \in [0,1]$.

\ENSURE $S \subset \mathbb{R}^n, \mathbf{w}=\mathbf{1}$.

\STATE \textbf{If}  $n=1$:  \\
Let $[a,b]$ be the smallest interval containing $(1-\eta-\epsilon)(1-\eta)$ fraction of the points,\\
$\tilde{S} \gets S\cap [a,b]$. \textbf{Return} $(\tilde{S},\mathbf{1})$

\STATE For each $d=1,\ldots, m$, 
\begin{enumerate}[(i)]
\item Let $[a_d,b_d]$ be the smallest interval containing $(1-\eta-\epsilon)(1-\eta)$ fraction of the points $\{x_d: \x \in S\},$
\item let $a_i \gets $MEAN$(\{x_d: \x \in S\} \cap [a_d,b_d])$. $\mathbf{a}\gets (a_1, \ldots, a_n)$.
\end{enumerate}

\STATE Set $B(r,\mathbf{a}) \gets$ ball of minium radius $r$ centered at $\mathbf{a}$ that contains $(1-\eta-\epsilon)(1-\eta)$ fraction of $S$.

\STATE $\tilde{S} \gets S\cap B(r,\mathbf{a})$. \textbf{Return} $(\tilde{S}, \mathbf{1})$.
\end{algorithmic}
\end{algorithm}

\begin{algorithm}[H]
\caption{AGNOSTICMEAN($\cdot$)}
\begin{algorithmic}[1]

\REQUIRE $S=\{ \x_1, \x_2,\ldots, \x_m\} \subset \mathbb{R}^n $.

\ENSURE $\hat{\mu}$.

\STATE Let $(\tilde{S}, \mathbf{w})$ = OUTLIERDAMPING(S) $ \sslash$ Gaussian Case \\
Let $(\tilde{S}, \mathbf{\omega})$ = OUTLIERREMOVAL(S,$1- w_1$) $ \sslash$ Non-Gaussian Case

\STATE \textbf{If}  $n=1$:  \\
\begin{enumerate}[(a)]
\item \textbf{if}  $\mathbf{\omega} = -\mathbf{1}$,  \textbf{Return} median$(\tilde{S})$.$ \sslash$ Gaussian Case 
\item \textbf{else} \textbf{Return} mean$(\tilde{S})$. $ \sslash$ Non-Gaussian Case
\end{enumerate}

\STATE $\mu_{\tilde{S}, \mathbf{w}} \gets \frac{1}{m} \sum_{i \in S} w_{x_i} \x_i$.\\
   $ \Sigma_{\tilde{S},\mathbf{w}} \gets \frac{1}{m} \sum_{i } w_{x_i}(\x_i - \mu_{\tilde{S}, \mathbf{w}})(\x_i - \mu_{\tilde{S}, \mathbf{w}})^T$.  \\
   Let $V$ be the span of the top $n/2$ principal components of $\Sigma_{\tilde{S},\mathbf{w}}$, and W be its complment.

\STATE  Set $S_1 \gets \mathbf{P}_V(S)$ where $\mathbf{P}_V$ is the projection operation on to $V$.

\STATE Set $\hat{\mu}_V \gets \text{AGNOSTICMEAN}(S_1)$ and $\hat{mu}_W \gets \text{MEAN}(\mathbf{P}_w \tilde{S})$

\STATE Let $\hat{\mu} \in \mathbb{R}^n$ be such that $\mathbf{P}_V \hat{\mu} = \hat{\mu}_V$, and $\mathbf{P}_W \hat{\mu} = \hat{\mu}_W$.

\STATE \textbf{Return} $\hat{\mu}$.
\end{algorithmic}
\end{algorithm}

\begin{algorithm}[H]
\caption{AGNOSTICCOV($S, w_2 + w_3$)}
\begin{algorithmic}[1]

\REQUIRE $S=\{ \x_1, \x_2,\ldots, \x_m\} \subset \mathbb{R}^n $.

\ENSURE $\hat{\mu}$.

\STATE $x'_i =\frac{x_i - x_{i+m/2}}{\sqrt{2}},$ for $i\in \{1,\ldots, m/2\}$. (Notice that $\E [x'(x')^T] =\Sigma$.)

\STATE Let $S^{(2)} \gets \{ x'_i (x_i')^T | i=1,2,\ldots, m/2\}$

\STATE  Run the mean estimation algorithm on $S^{(2)}$, where elements of $S^{(2)}$ are viewed as vectors in $\mathbb{R}^{n^2}$. Let the output be $\hat{\Sigma}$.

\STATE \textbf{Return} $\hat{\Sigma}$.

\end{algorithmic}
\end{algorithm}

\section{Proof of Lemma \ref{lemma1nons}}
\normalsize
Let $l(x)= (x-\hat{\mu}_1)' \hat{\Sigma}^{-1} (x-\hat{\mu}_1),$ and $${l^{\ast}}(x) = (x-\mu_1)'\Sigma^{-1} (x-\mu_1),$$ be the true value of the Mahalanobis distance. 
Then,
\begin{align*}
\begin{split}
l(x)  = &(x - \mu_1 + \mu_1 -\hat{\mu}_1)' \hat{\Sigma}_1^{-1} (x - \mu_1 + \mu_1 -\hat{\mu}_1) \\
  = & (x-\mu)'\hat{\Sigma}^{-1} (x-\mu_1) + 2 (\mu_1 -\hat{\mu}_1)' \hat{\Sigma}^{-1} (x - \mu_1) + ( \mu_1 -\hat{\mu}_1)' \hat{\Sigma}^{-1} (\mu_1 -\hat{\mu}_1) \\
  = & (x-\mu_1)'{\Sigma}^{-1} (x-\mu_1) + (x-\mu_1)' \hat{\Sigma}^{-1} (\Sigma_1 - \hat{\Sigma} ) \Sigma^{-1} (x-\mu_1) \\
 & + 2 (\mu_1 -\hat{\mu}_1)' \hat{\Sigma}^{-1} (x - \mu_1) + ( \mu_1 -\hat{\mu}_1)' \hat{\Sigma}^{-1} (\mu_1 -\hat{\mu}_1).
\end{split}
\end{align*}
And 
\begin{align*}
\begin{split}
& | (x-\mu_1)' \hat{\Sigma}^{-1} (\Sigma_1 - \hat{\Sigma} ) \Sigma^{-1} (x-\mu_1) | \\
\leq &   \| \Sigma - \hat{\Sigma} \|_F \| \hat{\Sigma}_1^{-1}\|_2  \|\Sigma^{-1}\|_2 \| x-\mu_1\|^2_2 \\
\leq & 2  \| \Sigma - \hat{\Sigma} \|_F \|\Sigma^{-1}\|^2_2 \| x-\mu_1\|^2_2 .
\end{split}
\end{align*}
 the last inequality is due to the fact that $\|X \|_2 \leq \|X \|_F$, and $\|(A+E)^{-1}\|_2 \leq \frac{\|A^{-1}\|_2}{1-\| A^{-1}E \|_2}$ as long as $A$ is non-singular and $\|A^{-1} E\|_2 <1$, such that
 \begin{align*}
 \begin{split}
 \| \hat{\Sigma}^{-1}\|_2  = &  \|( \Sigma + \hat{\Sigma} - \Sigma )^{-1} \| \\
 \leq &  \frac{\|\Sigma^{-1}\|_2}{ 1- \| \Sigma^{-1}(\Sigma - \hat{\Sigma}) \|_2} \\
 \leq &  \frac{\|\Sigma^{-1}\|_2}{ 1- \| \Sigma^{-1}\|_2 \| \Sigma- \hat{\Sigma} \|_F} \\
 \leq & 2 \|\Sigma^{-1}\|_2,
 \end{split}
 \end{align*}
{as long as} $ \| \Sigma - \hat{\Sigma} \|_F \leq \frac{1}{2} \| \Sigma^{-1}\|^{-1}_2$.
By Theorem 1.1, 1.3 from \cite{lai2016agnostic},
\begin{align*}
\begin{split}
\| \mu_1 - \hat{\mu}_1 \|_2 &= O(\eta^{1/2}+\epsilon)\|\Sigma\|_2^{1/2} \log^{1/2}n \\
\| \Sigma - \hat{\Sigma} \|_F &=O(\eta^{1/2} + C_{4,2}^{1/4}(\eta+\epsilon)^{3/4})C_4^{1/2}\|\Sigma\|_2 \log^{1/2}n,
\end{split}
\end{align*}
with $C_4= 3, C_{4,2} = 7!!/4$.
Therefore, 
\begin{align*}
\begin{split}
|l(x) - l^\ast(x)| \leq & 2  \| \Sigma^{-1}\|_2^2  \| x-\mu_1\|_2^2 \| \Sigma - \hat{\Sigma} \|_F + 4 \| \Sigma^{-1}\|_2 \|x-\mu_1\|_2  \| \mu_1 -\hat{\mu}_1\|_2 +  \| \Sigma^{-1}\|_2 \|\mu_1 -\hat{\mu}_1\|_2^2.
\end{split}
\end{align*}
Therefore, we can approximate $l(x)$ using $l^\ast(x)$ which follows (non-)central chi-squared distribution.


 Next, notice that
\begin{align}
\begin{split}
\label{bayesrule}
& \Pr( w \in G_1 | l(w) \geq l(X)_{(m_1)})  \\
 = & \frac{\Pr( l(w) \geq l(X)_{(m_1)} | w \in G_1 ) \Pr( w \in G_1) }{ \sum_{k=1}^2 \Pr(  l(w) \geq l(X)_{(m_1)} | w \in G_k) \Pr( w \in G_k)  +  \Pr(  l(w) \geq l(X)_{(m_1)} | w \in N) \Pr( w \in N) } \\
 \leq & \frac{\Pr( l(w) \geq l(X)_{(m_1)} | w \in G_1 ) \Pr( w \in G_1) }{  \Pr(  l(w) \geq l(X)_{(m_1)} | w \in G_1) \Pr( w \in G_1) + \Pr(  l(w) \geq l(X)_{(m_1)} | w \in G_2) \Pr( w \in G_2)  } .
  \end{split}
 \end{align}
 To show that the above is bounded, we need to find (i) an upper bound on $ \Pr( l(w) \geq l(X)_{(m_1)} | w \in G_1 )$,(ii) a lower bound for $\Pr(  l(w) \geq l(X)_{(m_1)} | w \in G_2)$ in \eqref{bayesrule}..
 
(i) First need to prove an upper bound on $ \Pr( l(w) \geq l(X)_{(m_1)} | w \in G_1 )$.
\begin{align*}
\begin{split}
 & \Pr( l(w) \geq l(X)_{(m_1)} | w \in G_1 )  \\
 \leq & \Pr( l(w) \geq l(X^{G_1\cup G_2})_{(m_1-m_3)} | w \in G_1  ) \\
  = & \Pr(l(w) \geq l(X^{G_1\cup G_2} )_{(m_1-m_3)}, l(w) \geq l(X^{G_1})_{(m_1(1-\beta) -m_3) } | w \in G_1  )  \\
  &+ \Pr(l(w) \geq l(X^{G_1\cup G_2} )_{(m_1-m_3)}, l(w) \leq l(X^{G_1})_{(m_1(1-\beta) -m_3 )} | w \in G_1  ) \\
  \leq & \beta +\frac{m_3}{m_1} + \sum_{n_1+n_2 = m_1 - m_3} \Pr( l(X^{G_1})_{(m_1(1-\beta) -m_3) } \geq  l(w)\geq l(X^{G_1})_{(n_1) }, l(w) \geq l(X^{G_2})_{(n_2)} | w \in G_1) \\
  \leq &  \beta +\frac{m_3}{m_1} + \sum_{n_1+n_2 = m_1 - m_3}\sum_{j=n_1}^{m_1(1-\beta)-m_3}  \Pr(  l(X^{G_1})_{(j) }=l(w) \geq l(X^{G_2})_{(n_2) } | w \in G_1) \\
  \leq  & \beta +\frac{m_3}{m_1} + \sum_{n_1+n_2 = m_1 - m_3}\sum_{j=n_1}^{m_1(1-\beta)-m_3}  \frac{1}{m_1} \Pr(  l(X^{G_1})_{(m_1(1-\beta) -m_3) } \geq l(X^{G_2})_{(m_1\beta) } | w \in G_1) \\
  \leq & \beta +\frac{m_3}{m_1} + \frac{(m_2 -m_1 \beta)^2}{2m_1} \Pr(  l(X^{G_1})_{(m_1(1-\beta) -m_3 )} \geq l(X^{G_2})_{(m_1\beta) }).
 \end{split}
 \end{align*}
 
 The last inequality above can be upper bounded as follows,
 \begin{align}
 \begin{split}
 \label{upperbound}
 &  \Pr(  l(X^{G_1})_{(m_1(1-\beta) -m_3 } \geq l(X^{G_2})_{(m_1\beta )}) \\
 \leq & \Pr(l(X^{G_1})_{(m_1(1-\beta) -m_3) }\geq l(X^{G_2})_{m_1\beta } | l(X^{G_2})_{(m_1\beta) } \geq t)\Pr(l(X^{G_2})_{m_1\beta } \geq t) \\
 & + \Pr( l(X^{G_2})_{(m_1\beta) } \leq t) \\
 \leq & \Pr(l(X^{G_1})_{(m_1(1-\beta) -m_3) }\geq t)\Pr(l(X^{G_2})_{(m_1\beta) } \geq t)  + \Pr( l(X^{G_2})_{(m_1\beta) } \leq t) \\
 \leq & \Pr( \frac{\sum_i \mathbbm{1}(l(X_i^{G_1}) \leq t)}{m_1} \leq 1-\beta-\frac{m_3}{m_1} ) + 
 \Pr(\frac{\sum_i \mathbbm{1}(l(X_i^{G_2}) \leq t)}{m_2} \geq \frac{m_1\beta}{m_2} ) \\
 \leq &\Pr(\frac{S_{m_1}(\underline{p})}{m_1} \leq 1-\beta-\frac{m_3}{m_1}) + \Pr(\frac{S_{m_2}(\bar{p})}{m_2} \geq \frac{m_1 \beta}{m_2}).
 \end{split}
 \end{align}
Here $\Pr(l(X_i^{G_1}) \leq t)  \geq \underline{p}$ and $\Pr(l(X_i^{G_2}) \leq t) \leq \bar{p}$ and $S_n(p)$ is a  sum of $n$ i.i.d. Bernoulli trials with $p$ being the probability of success. The last inequality is due to the following lemma. The proof for the lemma is standard calculus and thus omitted.
\begin{lemma}
Denote $S_m(p)$  the sum of $m$ Bernoulli random variables with parameter $p$. Then $\Pr(S_m(p) \leq t)$ is non-decreasing in $p$.
\end{lemma}

For $\underline{p}$,   let $\underline{p} = 1-\exp(-\frac{t'-\sqrt{(2t'-n)n}}{2})- \exp(-x')$, where 
\begin{align}
\begin{split}
\label{eq:t}
t'= & t-  2\|\Sigma -\hat{\Sigma}\|_F \|\Sigma^{-1}\|_2^2 (\tr(\Sigma) + 2\sqrt{\tr(\Sigma^2) x'} + 2 \|\Sigma\| x') \\
 & - 4 \|\mu_1 -\hat{\mu}_1\|_2 \|\Sigma^{-1}\|_2  \sqrt{\tr(\Sigma) + 2 \sqrt{\tr(\Sigma^2) x'} +2 \|\Sigma\| x'} \\
& - 2\|\mu_1 -\hat{\mu}_1\|_2^2 \|\Sigma^{-1}\|_2^2,
\end{split}
\end{align}
we have
 \begin{align*}
 \begin{split}
 \Pr( l(X_i^{G_1}) \leq t ) \geq &  \Pr(l^\ast (X_i^{G_1})\leq t') \\
\geq  & 1-\Pr( \chi^2_n \geq t' |  \| X_i^{G_1} -\mu_1\|_2  \leq \tr(\Sigma) + 2 \sqrt{\tr(\Sigma^2)x'} + 2\|\Sigma\|x'  ) - \exp(-x') \\
\geq & 1-\exp(-\frac{t'-\sqrt{(2t'-n)n}}{2}) - \exp(-x').
 \end{split}
 \end{align*}
 To get the second inequality above, we apply Proposition 1 in \cite{hsu2012tail} that for $Y\sim N(\mu,\Sigma)$,
 $$ \Pr( \| Y -\mu\|_2  \geq \tr(\Sigma) + 2 \sqrt{\tr(\Sigma^2)y} + 2\|\Sigma\|y) \leq \exp(-x).$$
The last inequality above is due to \cite{laurent2000adaptive} on the concentration inequality of  central chi-squared distribution. 
 \begin{lemma}[\cite{laurent2000adaptive}]
 \label{laurent2000}
 Let $X\sim \chi^2_n$, then for any $x\geq0$,
 $$ \Pr( X -n \geq 2\sqrt{xn} + 2x ) \leq \exp(-x),$$
 $$ \Pr( X -n \leq -2\sqrt{xn} ) \leq \exp(-x).$$
\end{lemma}
By simple calculation, we can deduce that $ \Pr( \chi^2_n \geq t) \leq \exp(- \frac{(\sqrt{2t-n} - \sqrt{n})^2}{4} ).$ Therefore, by the Chernoff-like bound in Lemma \ref{lemma: chernoff},  we have
$$ \Pr(\frac{S_{m_1}(\underline{p})}{m_1} \leq1-\beta-\frac{m_3}{m_1}) \leq \exp(-m_1 \frac{[\beta+\frac{m_3}{m_1} -\exp(-\frac{t'-\sqrt{(2t'-n)n}}{2}) - \exp(-x')]^2 }{2 (1-\beta -\frac{m_3}{m_1})(\beta + \frac{m_3}{m_1}) } ),$$
 when
  \begin{equation}
 \label{eq1}
 \beta+\frac{m_3}{m_1} \geq  2\exp(-\frac{t'-\sqrt{(2t'-n)n}}{2}) +2  \exp(-x').
 \end{equation}
 \begin{lemma}[Relative Entropy Chernoff Bound]
 \label{lemma: chernoff}
 Assume that $X\in[0,1]$, and $\E[X] = \mu$. Fix $\epsilon$, Define
 $$\text{MaxVar}[a,b] = \max_{p\in[a,b]} Var_p,$$
 where $Var_p$ is the variance of a random variable which is $1$ with probability $p$ and $0$ with probability $1 - p$. Then we have the following inequalities.
 $$\Pr(\bar{X}_n \geq \mu+\epsilon) \leq \exp( -n\frac{\epsilon^2}{2\text{MaxVar}[\mu,\mu+\epsilon]} ),$$
  $$\Pr(\bar{X}_n \leq \mu-\epsilon) \leq \exp( -n\frac{\epsilon^2}{2\text{MaxVar}[\mu-\epsilon,\mu]} ).$$
 \end{lemma}
The proof for the Relative Entropy Chernoff Bound 
is standard and is omitted here. Let $ x'= \log \frac{\beta + \frac{m_3}{m_1} }{4}$,  to show \eqref{eq1} it suffices to show that
 \begin{equation*}
 \label{eq2}
\frac{t'-\sqrt{(2t'-n)n}}{2}  \geq  \log(\beta/4 + \frac{m_3}{4m_1}),
 \end{equation*}
 which is equivalent to 
  \begin{equation}
 \label{eq3}
t' \geq \frac{\big(\sqrt{n}+ 2\sqrt{\log(\frac{\beta}{4} + \frac{m_3}{4m_1})}\big)^2 +n}{2}.
 \end{equation}

For the second term of \eqref{upperbound}, firstly, we show that $l^\ast(w)$ for any $w\in G_2$ follows a non-central chi-squared distribution. 

Suppose $\Sigma= Q^T D Q$ with $Q$ being real orthogonal matrix. Suppose that $Y\sim N(\mu_2,\Sigma)$,
\begin{align*}
(Y-\mu_1)^T\Sigma^{-1}(Y-\mu_1) &= (Y-\mu_1)^T Q^T D^{-1} Q(Y-\mu_1)  \\
&= [D^{-1/2}Q(Y-\mu_1)]^T [D^{-1/2}Q(Y -\mu_1)].
\end{align*}
Notice that $D^{-1/2}Q(Y-\mu_2)$ follows some Gaussian distribution with $\E[D^{-1/2}Q(Y-\mu_1)]=D^{-1/2}Q(\mu_2 - \mu_1)$, 
$$\text{Var}[D^{-1/2}Q(Y- \mu_1)] = D^{-1/2}Q (Q^T D Q) (D^{-1/2}Q)^T = I. $$
Hence, $D^{-1/2}Q(Y-\mu_2) \sim N(D^{-1/2}Q^T(\mu_2 -\mu_1), I)$. Therefore, $$l^\ast (w) = (w-\mu_2)^T\Sigma^{-1}(x-\mu_2) \sim \chi^2_n(\lambda),$$
where $\lambda= \| D^{-1/2}Q^T(\mu_2 -\mu_1)\|_2^2 = (\mu_2- \mu_1)^T \Sigma^{-1} (\mu_2- \mu_1)$.

Similarly we have, for some $t''$ to be specified later,
\begin{align*}
\begin{split}
\Pr(  l(X_i^{G_2}) \leq t ) & \leq \Pr(  \chi^2_n(\lambda) \leq t'' ) + \exp(-x'') \\
& \leq  \exp(- (\frac{n+\lambda - t''}{2\sqrt{n+2\lambda}})^2)  + \exp(-x'')\\
& = \bar{p},
\end{split}
\end{align*}
where \begin{align*}
\begin{split}
t''= & t +  2\|\Sigma -\hat{\Sigma}\|_F \|\Sigma^{-1}\|_2^2 (\tr(\Sigma) + 2 \sqrt{\tr(\Sigma^2) x''} + 2 \|\Sigma\| x'' +  \|\mu_2 - \mu_1\|^2 (  \frac{2\|\Sigma \| x''}{\sqrt{\tr(\Sigma^2)x''}} + 1) \\
 & + 4 \|\mu_1 -\hat{\mu}_1\|_2 \|\Sigma^{-1}\|_2  \sqrt{\tr(\Sigma) + 2 \sqrt{\tr(\Sigma^2) x''} + 2 \|\Sigma\| x +  \|\mu_2 - \mu_1\|^2 (  \frac{2\|\Sigma \| x''}{\sqrt{\tr(\Sigma^2)x''}} + 1)} \\
& + 2\|\mu_1 -\hat{\mu}_1\|_2^2 \|\Sigma^{-1}\|_2^2.
\end{split}
\end{align*}
The  result is due to  the following lemma.
\begin{lemma}[\cite{birge2001alternative}]
Let $X\sim \chi^2_n(\lambda)$, then for all $x>0$,
$$\Pr( X\geq (n+\lambda) + 2\sqrt{(n+2\lambda)x} + 2x) \leq \exp(-x),$$
$$\Pr( X\leq (n+\lambda) - 2\sqrt{(n+2\lambda)x} ) \leq \exp(-x).$$
\end{lemma}
 By Lemma \ref{lemma: chernoff}, $\Pr(\frac{S_{m_2}(\bar{p})}{m_2} \geq \frac{m_1 \beta}{m_2}) \leq  \exp( -m_2 \frac{[\beta \frac{m_1}{m_2} -  \exp(- (\frac{n+\lambda - t''}{2\sqrt{n+2\lambda}})^2) - \exp(-x'') ]^2}{ \beta \frac{m_1}{m_2} ( 1-\beta \frac{m_1}{m_2})  } ) $ when
\begin{equation}
\label{eq4}
\beta \frac{m_1}{m_2} \geq 2  \exp(- (\frac{n+\lambda - t''}{2\sqrt{n+2\lambda}})^2)  +  2\exp(-x'').
\end{equation}
where $x'' =- \log \frac{\eta}{4}$. Then it suffices to show
$\beta \frac{m_1}{m_2} \geq 4  \exp(- (\frac{n+\lambda - t''}{2\sqrt{n+2\lambda}})^2),$
which is equivalent to,
\begin{equation}
\label{eq6}
t'' \leq n+\lambda - 2\sqrt{(n+2\lambda)\log\frac{m_2}{m_1\beta}}.
\end{equation}

Then we have 
 $$  \Pr(  l(X^{G_1})_{m_1(1-\beta) -m_3 } \geq l(X^{G_2})_{m_1\beta }) \sim O(\exp(-m_1 \beta)). $$ 
 Therefore, let $\beta = \frac{m_2}{m_1}\eta ,$ we have
\begin{align*}
\begin{split}
& \beta +\frac{m_3}{m_1} + \frac{(m_2 -m_1 \beta)^2}{2m_1} \Pr(  l(X^{G_1})_{m_1(1-\beta) -m_3 } \geq l(X^{G_2})_{m_1\beta })  \\
= &  \frac{m_2}{m_1} \eta + \frac{m_3}{m_1}  + c \frac{m_2}{m_1} \exp ( -m_2 \eta + \log m_2).
\end{split}
\end{align*}
as long as $\eta > \frac{\log m_2}{m_2} $, we have 
$$\frac{m_2}{m_1} \eta + \frac{m_3}{m_1}  + c \frac{m_2}{m_1} \exp ( -m_2 \eta + \log m_2)= o(\frac{m_2}{m_1}).$$

The only thing left is to find $t',t''$ that satisfy \eqref{eq3} and \eqref{eq6}. By  \eqref{eq3} 
\begin{align}
\begin{split}
\label{eq7}
t \geq  & 2\|\Sigma -\hat{\Sigma}\|_F \|\Sigma^{-1}\|_2^2 (\tr(\Sigma) + 2\sqrt{- \tr(\Sigma^2)  \log (\frac{\beta}{4} + \frac{m_3}{4m_1})} - 2 \|\Sigma\|  \log (\frac{\beta}{4} + \frac{m_3}{4m_1}) ) \\
 & + 4 \|\mu_1 -\hat{\mu}_1\|_2 \|\Sigma^{-1}\|_2  \sqrt{-\tr(\Sigma) - 2\sqrt{\tr(\Sigma^2)  \log (\frac{\beta}{4} + \frac{m_3}{4m_1})} - 2 \|\Sigma\|  \log (\frac{\beta}{4} + \frac{m_3}{4m_1})} \\
& + 2\|\mu_1 -\hat{\mu}_1\|_2^2 \|\Sigma^{-1}\|_2^2\\
& + \frac{\big(\sqrt{n}+ \sqrt{-2\log(\frac{\beta}{4} + \frac{m_3}{4m_1})}\big)^2 +n}{2}.
\end{split}
\end{align}

For \eqref{eq6}, we have
\begin{align*}
\begin{split}
& n+\lambda - 2\sqrt{(n+2\lambda)\log\frac{m_2}{m_1\beta}}  \\
\geq & t +  2\|\Sigma_1 -\hat{\Sigma}_1\|_F \|\Sigma_1^{-1}\|_2^2 (\tr(\Sigma_2) + 2 \sqrt{\tr(\Sigma_2^2)  \log \frac{4}{\eta}} + 2 \|\Sigma_2\| \log \frac{4}{\eta} +  \|\mu_2 - \mu_1\|^2 (  \frac{2\|\Sigma_2 \| \sqrt{ \log \frac{4}{\eta}} }{\sqrt{\tr(\Sigma_2^2)}} + 1) \\
 & + 4 \|\mu_1 -\hat{\mu}_1\|_2 \|\Sigma_1^{-1}\|_2  \sqrt{\tr(\Sigma_2) + 2 \sqrt{\tr(\Sigma_2^2) \log \frac{4}{\eta}} + 2 \|\Sigma_2\| \log \frac{4}{\eta} +  \|\mu_2 - \mu_1\|^2 (  \frac{2\|\Sigma_2 \| \sqrt{\log \frac{4}{\eta}}}{\sqrt{\tr(\Sigma_2^2)}} + 1)} \\
& + 2\|\mu_1 -\hat{\mu}_1\|_2^2 \|\Sigma_1^{-1}\|_2^2,
\end{split}
\end{align*}
combined with the lower bound on $t$ in \eqref{eq7}, we have 
\begin{align}
\label{eq8}
\begin{split}
& n+\lambda - 2\sqrt{(n+2\lambda)\log\frac{1}{\eta}} \\
\geq & 2\|\Sigma -\hat{\Sigma}\|_F \|\Sigma^{-1}\|_2^2 \times \big [ 2\tr(\Sigma) + 2 \sqrt{\tr(\Sigma^2) }(\sqrt{ \log \frac{4}{\eta}} + \sqrt{-\log (\frac{\eta}{4}\frac{m_2}{m_1} + \frac{m_3}{4m_1})} ) \\
& + 2 \|\Sigma\|( \log \frac{4}{\eta}- \log (\frac{\eta}{4}\frac{m_2}{m_1} + \frac{m_3}{4m_1})) +  \|\mu_2 - \mu_1\|^2 (  \frac{2\|\Sigma \| \sqrt{ \log \frac{4}{\eta}} }{\sqrt{\tr(\Sigma^2)}} + 1) ] \\
 & + 4 \|\mu_1 -\hat{\mu}_1\|_2 \|\Sigma^{-1}\|_2  \sqrt{\tr(\Sigma) + 2 \sqrt{\tr(\Sigma^2) \log \frac{4}{\eta}} + 2 \|\Sigma\| \log \frac{4}{\eta} +  \|\mu_2 - \mu_1\|^2 (  \frac{2\|\Sigma \| \sqrt{\log \frac{4}{\eta}}}{\sqrt{\tr(\Sigma^2)}} + 1)} \\
  & + 4 \|\mu_1 -\hat{\mu}_1\|_2 \|\Sigma^{-1}\|_2  \sqrt{\tr(\Sigma) + 2\sqrt{-\tr(\Sigma^2)  \log (\frac{\eta}{4} \frac{m_2}{m_1} + \frac{m_3}{4m_1})} - 2 \|\Sigma\|  \log (\frac{\eta}{4}\frac{m_2}{m_1} + \frac{m_3}{4m_1})} \\
& + 4 \|\mu_1 -\hat{\mu}_1\|_2^2 \|\Sigma^{-1}\|_2^2  + \frac{\big(\sqrt{n}+ \sqrt{2\log(\frac{\beta}{4} + \frac{m_3}{4m_1})}\big)^2 +n}{2}.
\end{split}
\end{align}

Notice that
$$\log\frac{4}{\eta} = O(- \log(\frac{\eta}{4}\frac{m_2}{m_1} + \frac{m_3}{4m_1})),$$
and
$$\lambda= (\mu_2 - \mu_1)' \Sigma^{-1}  (\mu_2 - \mu_1) \geq \|\Sigma\|_2^{-1} \|\mu_2 - \mu_1\|^2_2.$$
Hence, to show \eqref{eq8} it suffices to show that
\begin{align}
\label{eq9}
\begin{split}
\small
& n+\lambda - 2\sqrt{(n+2\lambda)\log\frac{1}{\eta}} - 2\|\Sigma -\hat{\Sigma}\|_F \|\Sigma^{-1}\|_2^2 \lambda (  \frac{2\|\Sigma \|^2 \sqrt{ \log \frac{4}{\eta}} }{\sqrt{\tr(\Sigma^2)}} + 1) \\
&  -8 \|\mu_1 -\hat{\mu}_1\|_2 \|\Sigma^{-1}\|_2  \sqrt{\lambda (  \frac{2\|\Sigma \|^2 \sqrt{\log \frac{4}{\eta}}}{\sqrt{\tr(\Sigma^2)}} + 1)} \\
\geq & 2\|\Sigma -\hat{\Sigma}\|_F \|\Sigma^{-1}\|_2^2 \times \big [ 2\tr(\Sigma) + 4 \sqrt{\tr(\Sigma^2) } \sqrt{\log (\frac{\eta}{4}\frac{m_2}{m_1} + \frac{m_3}{4m_1})^{-1}} ) + 4 \|\Sigma\| \log (\frac{\eta}{4}\frac{m_2}{m_1} + \frac{m_3}{4m_1})^{-1} ] \\
 & + 8 \|\mu_1 -\hat{\mu}_1\|_2 \|\Sigma^{-1}\|_2  \sqrt{2\tr(\Sigma) + 4 \sqrt{\tr(\Sigma^2) \log (\frac{\eta}{4} \frac{m_2}{m_1} + \frac{m_3}{4m_1})^{-1}} + 4 \|\Sigma\|\log (\frac{\eta}{4} \frac{m_2}{m_1} + \frac{m_3}{4m_1})^{-1} }\\
& + 4 \|\mu_1 -\hat{\mu}_1\|_2^2 \|\Sigma^{-1}\|_2^2  + \frac{\big(\sqrt{n}+ \sqrt{2\log(\frac{\eta}{4}\frac{m_2}{m_1} + \frac{m_3}{4m_1})^{-1}}\big)^2 +n}{2}.
\end{split}
\end{align}


Subsequently, \eqref{eq9} can be reduced to the following lower bound on $\lambda$.
\begin{align}
\begin{split}
\label{eq10}
\sqrt{\lambda} \geq & a + \sqrt{a^2 + 2b}= O(\sqrt{a^2 + 2b}), \\
a= &  8 \|\mu_1 -\hat{\mu}_1\|_2 \|\Sigma^{-1}\|_2  \sqrt{   \frac{2\|\Sigma \|^2 \sqrt{\log \frac{4}{\eta}}}{\sqrt{\tr(\Sigma^2)}} + 1},\\
b= & 2\|\Sigma -\hat{\Sigma}\|_F \|\Sigma^{-1}\|_2^2 \times \big [ 2\tr(\Sigma) + 4 \sqrt{\tr(\Sigma^2) } \sqrt{\log (\frac{\eta}{4}\frac{m_2}{m_1} + \frac{m_3}{4m_1})^{-1}} ) + 4 \|\Sigma\| \log (\frac{\eta}{4}\frac{m_2}{m_1} + \frac{m_3}{4m_1})^{-1} ] \\
 & + 8 \|\mu_1 -\hat{\mu}_1\|_2 \|\Sigma^{-1}\|_2  \sqrt{2\tr(\Sigma) + 4 \sqrt{\tr(\Sigma^2) \log (\frac{\eta}{4} \frac{m_2}{m_1} + \frac{m_3}{4m_1})^{-1}} + 4 \|\Sigma\|\log (\frac{\eta}{4} \frac{m_2}{m_1} + \frac{m_3}{4m_1})^{-1} }\\
& + 4 \|\mu_1 -\hat{\mu}_1\|_2^2 \|\Sigma^{-1}\|_2^2  + \sqrt{2n\log(\frac{\eta}{4}\frac{m_2}{m_1} + \frac{m_3}{4m_1})^{-1}} + \log(\frac{\eta}{4}\frac{m_2}{m_1} + \frac{m_3}{4m_1})^{-1}.
\end{split}
\end{align}

To reduce this complicated inequality,  it can be shown that,
$$a  \leq O( (w_2 \log n)^{\frac{1}{4}}  ).$$
Combined with
$$\|\Sigma -\hat{\Sigma}\|_F = O(\sqrt{w_2}) \|\Sigma\| \sqrt{\log n}, \quad \|\mu_1 -\hat{\mu}_1\|_F = O(\sqrt{w_2}) \|\Sigma\|^{1/2} \sqrt{\log n},$$
we have
$$ b= c_1 \sqrt{w_2\log n} ( \|\Sigma\|^{\frac{1}{2}} \|\Sigma^{-1}\| +1) ( \tr(\Sigma)  +  \sqrt{\tr(\Sigma^2)\log\frac{1}{\delta}} + \|\Sigma\| \log\frac{1}{\delta} ) + \frac{1}{\sqrt{2}} \sqrt{n\log \frac{1}{\delta}}  + \log\frac{1}{\delta},$$ 
where $\eta \frac{w_2}{w_1} + \frac{w_2}{w_3}=\delta,$ and  $c_1>0$ being some constant.
\normalsize Hence, by \eqref{eq10},
$$\lambda \geq c_1 \sqrt{w_2\log n} ( \|\Sigma\|^{\frac{1}{2}} \|\Sigma^{-1}\| +1) ( \tr(\Sigma)  +  \sqrt{\tr(\Sigma^2)\log\frac{1}{\delta}} + \|\Sigma\| \log\frac{1}{\delta})  + \frac{1}{\sqrt{2}} \sqrt{n\log \frac{1}{\delta}}  + \log\frac{1}{\delta}.$$

\normalsize
If $n=o(\log \frac{1}{\epsilon})$, then we have 
\begin{align}
\begin{split}
\lambda \geq \big(1+ c \sqrt{w_2 \log n}( \|\Sigma\|^{\frac{1}{2}} \|\Sigma^{-1}\| +1) \|\Sigma\| \big) \log\frac{1}{\epsilon} .
\end{split}
\end{align}

With this lower bound on $\lambda$, we can choose $t',t''$ to complete the upper bound of \eqref{upperbound}.

(ii) Next, we show a lower bound for $\Pr(  l(w) \geq l(X)_{(m_1)} | w \in G_2)$ in \eqref{bayesrule}. 
\begin{align*}
\begin{split}
& \Pr( l(w) \geq l(X)_{(m_1)} | w\in G_2) \\
= &   \sum_{i + j+k= m_1} \Pr(l(w) \geq l(X^{G_1})_{(i)}, l(w) \geq l(X^{G_2})_{(j)}, l(w) \geq l(X^{N})_{(k)} | w\in G_2)\\
\geq &  \sum_{i + j= m_1} \Pr(l(w) \geq l(X^{G_1})_{(i)}, l(w) \geq l(X^{G_2})_{(j)} | w\in G_2) \\
 \geq & \Pr(l(w) > l(X^{G_1})_{(m_1 - \tilde{m}_2 )}, l(w) \geq l(X^{G_2})_{(\tilde{m}_2 )} | w\in G_2, l(X^{G_1})_{(m_1 - \tilde{m}_2 )} \leq t )  \\
 & \times \Pr(  l(X^{G_1})_{(m_1 - \tilde{m}_2 )} \leq t ) \\
 \geq & \Pr(l(w) > t , l(w) \geq l(X^{G_2})_{(\tilde{m}_2)} | w\in G_2 ) \times \Pr(  l(X^{G_1})_{(m_1 - \tilde{m}_2 )} \leq t ) \\ 
 = & \frac{m_2 - \tilde{m}_2}{m_2} \Pr(l(X^{G_2})_{(\tilde{m}_2 )} > t ) \times \Pr(  l(X^{G_1})_{(m_1 - \tilde{m}_2 )} \leq t ).
\end{split}
\end{align*}
for some $t, \tilde{m}_2$ to be specified later. 

\noindent\textbf{Claim 1.}  $\Pr(  l(X^{G_1})_{(m_1 - \tilde{m}_2 )} \geq t ) \leq  O(\exp [ -\tilde{m}_2 (1-\frac{\bar{p}}{ \tilde{m}_2/m_1})^2 ] )    $. \\
\textit{Proof.} 
\begin{align*}
\begin{split}
\Pr(l(X^{G_1})_{(m_1 - \tilde{m}_2 )} \geq t) & = \Pr(\sum_{i=1}^{m_1} \mathbbm{1}(l(X_i^{G_1}) \geq t) \geq \tilde{m}_2) \\
& =  \Pr(\frac{ \sum_{i=1}^{m_1} \mathbbm{1}(l(X_i^{G_1}) \geq t)}{m_1} \geq \frac{\tilde{m}_2}{m_1} ) \\
& = \Pr(\frac{S_{m_1}(p)}{m_1} \geq \frac{\tilde{m}_2}{m_1} ) \\
& \leq \Pr(\frac{S_{m_1}(\bar{p})}{m_1} \geq \frac{\tilde{m}_2}{m_1}) \\
& \leq  \exp[-m_1 \frac{ (\frac{\tilde{m}_2}{m_1}- \bar{p})^2 }{2 \text{MaxVar}[\bar{p},\frac{\tilde{m}_2}{m_1}] }]. 
 \end{split}
\end{align*}
By Lemma \ref{laurent2000}, $ \Pr( \chi^2_n \geq t) \leq \exp(- \frac{(\sqrt{2t-n} - \sqrt{n})^2}{4} ).$
Let $\bar{p} =\exp(- \frac{(\sqrt{2t-n} - \sqrt{n})^2}{4} ) =  \exp(- \frac{t - \sqrt{(2t-n)n}}{2})$, we have 
\begin{align*}
\begin{split}
\Pr(l(X^{G_1})_{(m_1 - \tilde{m}_2 )} \geq t) & \leq  \exp[-m_1 \frac{ (\frac{\tilde{m}_2}{m_1}- \bar{p})^2 }{2 \text{MaxVar}[\bar{p},\frac{\tilde{m}_2}{m_1}] }] \\
& = \exp[-m_1 \frac{ (\frac{\tilde{m}_2}{m_1}- \bar{p})^2 }{2 \frac{\tilde{m}_2}{m_1}(1- \frac{\tilde{m}_2}{m_1})  ] }] \\
&= O(\exp [ -\tilde{m}_2 (1-\frac{\bar{p}}{ \tilde{m}_2/m_1})^2 ] ).           
 \end{split}
\end{align*}
under 
\begin{equation}
\label{eq:claim1ineq}
\frac{\tilde{m}_2}{m_1} > \bar{p} = \exp(- \frac{t - \sqrt{(2t-n)n}}{2}).
\end{equation}

\noindent\textbf{Claim 2.} $ \Pr(l(X^{G_2})_{(\tilde{m}_2 )} > t ) \geq 1/2$. \\
\textit{Proof.} 
\begin{align*}
\begin{split}
\Pr ( l(X^{G_2})_{(\tilde{m}_2 )}) \leq t) = &  \Pr(\sum_{i=1}^{m_2} \mathbbm{1}(l(X_i^{G_2}) \leq t) \geq  \tilde{m}_2) \\
= & \Pr( \frac{S_{m_2}(q)}{m_2} \geq \frac{\tilde{m}_2}{m_2} ) \\
\leq & \Pr( \frac{S_{m_2}(\bar{q})}{m_2} \geq \frac{\tilde{m}_2}{m_2} ). 
\end{split}
\end{align*}
Need to find $\tilde{m}_2, \bar{q}$ such that 
$$\frac{\tilde{m}_2}{m_2} \geq \bar{q} \geq \Pr(l(X^{G_2}) \leq t)), $$ which then leads to $\Pr( \frac{S_{m_2}(\bar{q})}{m_2} \geq \frac{\tilde{m}_2}{m_2} )  \leq \frac{1}{2}$.
According to  \cite{birge2001alternative}, $\forall t < n+\lambda$,
$$\Pr( \chi_n^2(\lambda) \leq t) \leq \exp( - ( \frac{n+\lambda - t}{2\sqrt{n+2\lambda} } )^2) = \bar{q} . $$
Therefore, combined with \eqref{eq:claim1ineq} in Claim 1,  we have 
$$\frac{\tilde{m}_2}{m_2} \geq  \max \{ \exp( - ( \frac{n+\lambda - t}{2\sqrt{n+2\lambda} } )^2),\exp(- \frac{t - \sqrt{(2t-n)n}}{2} + \log \frac{m_1}{2m_2}) \}, $$
which is equivalent to, 
$$ n+ \lambda - 2\sqrt{ -\log \frac{\tilde{m}_2}{m_2}  }  \sqrt{n+ 2 \lambda} \geq t \geq \frac{(2\sqrt{\log \frac{m_2}{\tilde{m}_2}+  \log \frac{m_1}{2m_2}} +\sqrt{n})^2 + n }{2}. $$
Hence,  if
\begin{align*}
\begin{split}
n+ \lambda - 2\sqrt{ -\log \frac{\tilde{m}_2}{m_2}  }  \sqrt{n+ 2 \lambda} & \geq \frac{(2\sqrt{\log \frac{m_2}{\tilde{m}_2}+  \log \frac{m_1}{2m_2}} +\sqrt{n})^2 + n}{2}, 
\end{split}
\end{align*}
i.e.,  
\begin{equation}
\lambda \geq 2 (\sqrt{\log \frac{m_2}{\tilde{m}_2}+  \log \frac{m_1}{2m_2}} + \sqrt{ \log \frac{m_2}{\tilde{m}_2}}) ^2 + 2 ( \sqrt{\log \frac{m_2}{\tilde{m}_2}+  \log \frac{m_1}{2m_2}} + \sqrt {\log {\frac{m_2}{\tilde{m}_2}}} ) \sqrt{n}.
\end{equation}
Let $\tilde{m}_2 = \frac{1}{2} m_2$ we have,
\begin{equation}
\label{eq:separationclaim1}
\lambda \geq 3\log \frac{m_1}{2m_2}  + 3 \sqrt{n\log \frac{m_1}{2m_2} },
\end{equation}
which suffices to prove Claim 2.

The following presents an alternative way to prove Claim 2. In particular, it relies on the following characterization of a non-central chi-squared distribution.

\begin{lemma}
$$\chi_n^2(\lambda) \myeq \chi_{n-1}^2(0) +  \chi_1^2(\lambda).$$
\end{lemma}
The CDF of non-central chi-square with one degree of freedom is quite explicit.
$$ \Pr( \chi^2_1(\lambda) \leq t) = \Phi(\sqrt{t} - \sqrt{\lambda}) - \Phi(-\sqrt{t} - \sqrt{\lambda}). $$
Hence,
\begin{align*}
\begin{split}
\Pr ( l(X^{G_2})_{(\tilde{m}_2 )}) \leq t) & = \Pr(\sum_{i=1}^{m_2} \mathbbm{1}(l(X_i^{G_2}) \geq t) \leq m_2 -  \tilde{m}_2) \\
& = \Pr( \frac{S_{m_2}(p) }{m_2} \leq \frac{m_2 - \tilde{m}_2}{m_2}  ) \\
& \leq  \Pr( \frac{S_{m_2}(\underline{p}) }{m_2} \leq \frac{m_2 - \tilde{m}_2}{m_2}  ),
\end{split}
\end{align*}
where $p = \Pr( \chi^2_n (\lambda) \geq t) \geq \underline{p}. $
Need to find a lower bound $\underline{p}$ and $\tilde{m}_2$, such that $\underline{p} \geq \frac{m_2 - \tilde{m}_2}{m_2}$ and thus
$$\Pr ( l(X^{G_2})_{(\tilde{m}_2)}) \leq t) \leq \frac{1}{2},$$
which then completes the proof of Claim 2.

 Denote by  $\bar{\Phi}(t) = 1-\Phi(t)$. By Lemma 3.6 and the explicit CDF of non-central chi-square with one degree of freedom,
\begin{align*}
\begin{split}
& \Pr(\chi^2_1(\lambda) + \chi^2_{n-1}(0) \geq t) \\
\geq & \Pr(\chi^2_1(\lambda) + \chi^2_{n-1}(0) \geq t | \chi^2_{n-1}(0) \geq {n-1} - 2\sqrt{(n-1)x}) \Pr(\chi^2_{n-1}(0) \geq {n-1} - 2\sqrt{(n-1)x}) \\
\geq & \Pr(\chi^2_1(\lambda) \geq  t- (n-1) + 2\sqrt{(n-1)x}) \Pr(\chi^2_{n-1}(0) \geq {n-1} - 2\sqrt{(n-1)x})  \\
\geq & \bar{\Phi}( \sqrt{t- (n-1) + 2\sqrt{(n-1)x}} - \sqrt{\lambda} ) (1-\exp(-x)) \\
= & \underline{p} .
\end{split}
\end{align*}
Combined with \eqref{eq:claim1ineq} in Claim 1, it suffices to find $t, \tilde{m}_2$  satisfying
\begin{align}
\begin{split}
\label{eq:m2claim2}
 \frac{\tilde{m}_2}{m_2} \geq &1-\underline{p}, \\
\frac{\tilde{m}_2}{m_2} \geq  &\exp(- \frac{t - \sqrt{(2t-n)n}}{2} + \log \frac{m_1}{2m_2}).
\end{split}
\end{align}
Let $t- (n-1) + 2\sqrt{(n-1)x }= \lambda + 1$, $x=1$, and 
$$\frac{3}{4} \geq \frac{\tilde{m}_2}{m_2} = 1- \bar{\Phi}(1)(1-e^{-1}) \geq \frac{1}{2},$$ which yields $  \Pr(\chi^2_1(\lambda) + \chi^2_{n-1}(0) \geq t) \geq \underline{p} = \bar{\Phi}(1)(1-e^{-1})$. Moreover, we need to show that
\begin{align*}
\begin{split}
\frac{t - \sqrt{(2t-n)n}}{2} - \log \frac{m_1}{2m_2} & \geq - \log \frac{\tilde{m}_2}{m_2} =  \log (1-\underline{p})^{-1}.
\end{split}
\end{align*}
by plugging in $t- (n-1) + 2\sqrt{(n-1)}= \lambda + 1$, 
$\lambda \geq \frac{(\sqrt{n}+ 2\sqrt{\log\frac{m_1}{m_2}})^2 +n}{2}- n + 2\sqrt{n-1}$, which finally yields,
$$\lambda \geq 2 \sqrt{n\log\frac{m_1}{m_2}} + 2\sqrt{n-1} + 2\log\frac{m_1}{m_2}.$$
Notice that this is actually quite close to \eqref{eq:separationclaim1}. Hence, all of the above combined guarantees that $\tilde{m}_2$ exists, 
and thus completes the alternative proof of Claim 2.

Therefore,
\begin{align*}
\begin{split}
 \Pr( l(w) \geq l(X)_{(m_1)} | w\in G_2)\geq & \frac{m_2 - \tilde{m}_2}{m_2} \Pr(l(X^{G_2})_{(m_2+1-\tilde{m}_2 : m_2)} > t ) \times \Pr(  l(X^{G_1})_{(m_1 - \tilde{m}_2 )} \leq t ) \\
 \geq & \frac{1}{20}(1-\exp(-c m_2)),
\end{split}
\end{align*}
for some $c>0$.

Moreover, given enough sampled points, the fraction of sampled points belonging to  each component is approximately the true weight $w_i$, that is, given
$m \geq \frac{2 \log{\frac{1}{\eta}}} {w_2 \eta^2 }$, we have $ \exp(-\frac{m w_2 \eta^2 } {2})\leq \eta,$
which yields the concentration inequality of Bernoulli trials,
$$\Pr(| \frac{m_2}{m} - w_2 | \geq w_2 {\eta}) \leq \exp(- m w_2 \eta^2 ) \leq \eta. $$
Similarly,
with $m \geq \frac{2 \log{\frac{1}{\eta}}} {w_3}$, we can bound the probability of Bernoulli trials,
$$\Pr(| \frac{m_3}{m} - w_3 | \geq  w_3 ) \leq 2\exp(-2 m w_3 ).$$
Hence, by union bounds,
\begin{align*}
\begin{split}
& \Pr( w \in G_1 | l(w) \geq l(X)_{(m_1)})  \\
 = & \frac{\Pr( l(w) \geq l(X)_{(m_1)} | w \in G_1 ) \Pr( w \in G_1) }{ \sum_{k=1}^2 \Pr(  l(w) \geq l(X)_{(m_1)} | w \in G_k) \Pr( w \in G_k)  +  \Pr(  l(w) \geq l(X)_{(m_1)} | w \in N) \Pr( w \in N) } \\
 \leq & \frac{\Pr( l(w) \geq l(X)_{(m_1)} | w \in G_1 ) \Pr( w \in G_1) }{  \Pr(  l(w) \geq l(X)_{(m_1)} | w \in G_1) \Pr( w \in G_1) + \Pr(  l(w) \geq l(X)_{(m_1)} | w \in G_2) \Pr( w \in G_2)  }   \\
 \leq & \frac{ 2\eta w_2 w_1 }{ 2\eta  w_2 w_1 + w_2  \frac{1}{20}} + 2 \exp(-\frac{m (w_2 \eta)^2 } {2w_2}) + 2 \exp(-\frac{m w_3 } {2})  \\
 \leq & 32 \eta. 
  \end{split}
 \end{align*}
To conclude, under the condition of Lemma \ref{lemma1nons} and Lemma 3.2, 3.3, in order to complete the proof of Theorem \ref{thm1}, $m$ must satisfy,
$$m = \Omega (\frac{n(\frac{1}{w_2} +n)( \log n + \log \frac{1}{\epsilon}) \log n}{\epsilon^2} ).$$
$$m \geq \frac{2 \log \frac{1}{\epsilon^2 + \frac{w_3}{w_2} } }{w_2 ( \epsilon^2 + \frac{w_3}{w_2})^2} .$$
$$ m \geq \frac{2 \log \frac{1}{\epsilon^2 + \frac{w_2}{w_3} } }{w_3}. $$
Hence, 
$$ m = \Omega ( \frac{n(n+\frac{1}{w_2})( \log n + \log \frac{1}{\epsilon}) \log n}{\epsilon^2} +  \frac{ \frac{1}{w_2} \log (\frac{1}{\epsilon} + \sqrt{\frac{w_2}{w_3} }) }{\epsilon^4 + (\frac{w_3}{w_2})^2 } +  \frac{ \log (\frac{1}{\epsilon} + \sqrt{\frac{w_2}{w_3} }) }{w_3} ) .$$

 \section{Proof of Lemma \ref{lemma1s}}
 Let ${l^{\ast}}(x) = \frac{\|x-\mu_1\|^2_2}{\sigma_1^2}$  be the true value of Mahalanobis distance, and $l(x)= \frac{\|x-\hat{\mu}_1\|_2^2}{\hat{\sigma}^2}$ the computed estimator for $l(x)$ in the algorithm.
Then 
\begin{align*}
\begin{split}
l(x)   = & \frac{\|x-\mu_1\|^2_2}{\hat{\sigma}_1^2} + 2 \frac{(\mu_1 - \hat{\mu}_1)'(x-\mu_1)}{{\hat{\sigma}_1^2}} + \frac{\|\mu_1 - \hat{\mu}_1\|_2^2}{\hat{\sigma}_1^2}.
\end{split}
\end{align*}
$ \hat{\sigma}^2_1 \leq 2  {\sigma}_1^2$
By Theorem 1.1, 1.3 from \cite{lai2016agnostic},
\begin{align*}
\begin{split}
\| \mu_1 - \hat{\mu}_1 \|_2 &= O(\eta+ \epsilon) \sigma_1\sqrt{ \log n}. \\
| \sigma^2_1 - \hat{\sigma}^2_1 | &=O(\eta+\epsilon) \sigma_1^2.
\end{split}
\end{align*}
Similarly by Bayes Rule,
\begin{align*}
\begin{split}
& \Pr( w \in G_1 | l(w) \geq l(X)_{(m_1)})  \\
 = & \frac{\Pr( l(w) \geq l(X)_{(m_1)} | w \in G_1 ) \Pr( w \in G_1) }{ \sum_{k=1}^2 \Pr(  l(w) \geq l(X)_{(m_1)} | w \in G_k) \Pr( w \in G_k)  +  \Pr(  l(w) \geq l(X)_{(m_1)} | w \in N) \Pr( w \in N) } \\
 \leq & \frac{\Pr( l(w) \geq l(X)_{(m_1)} | w \in G_1 ) \Pr( w \in G_1) }{  \Pr(  l(w) \geq l(X)_{(m_1)} | w \in G_1) \Pr( w \in G_1) + \Pr(  l(w) \geq l(X)_{(m_1)} | w \in G_2) \Pr( w \in G_2)  } .
  \end{split}
 \end{align*}
Need to prove an upper bound on $ \Pr( l(w) \geq l(X)_{(m_1)} | w \in G_1 )$ and a lower bound on $\Pr(  l(w) \geq l(X)_{(m_1)} | w \in G_2)$. Similar to the non-spherical case,
\begin{align*}
\begin{split}
  \Pr( l(w) \geq l(X)_{(m_1)} | w \in G_1 )  \leq \beta +\frac{m_3}{m_1} + \frac{(m_2 -m_1 \beta)^2}{2m_1} \Pr(  l(X^{G_1})_{(m_1(1-\beta) -m_3)} \geq l(X^{G_2})_{(m_1\beta)}),
 \end{split}
 \end{align*}
and we have
 \begin{align*}
 \begin{split}
 &  \Pr(  l(X^{G_1})_{(m_1(1-\beta) -m_3) } \geq l(X^{G_2})_{(m_1\beta)}) \\
 \leq & \Pr(l(X^{G_1})_{m_1(1-\beta) -m_3 :m_1}\geq l(X^{G_2})_{m_1\beta :m_2} | l(X^{G_2})_{m_1\beta :m_2} \geq t)\Pr(l(X^{G_2})_{m_1\beta :m_2} \geq t) \\
 & + \Pr( l(X^{G_2})_{m_1\beta :m_2} \leq t) \\
  \leq &\Pr(\frac{S_{m_1}(\underline{p})}{m_1} \leq 1-\beta-\frac{m_3}{m_1}) + \Pr(\frac{S_{m_2}(\bar{p})}{m_2} \geq \frac{m_1 \beta}{m_2}).
 \end{split}
 \end{align*}
Here $\Pr(l(X_i^{G_1}) \leq t)  \geq \underline{p}$ and $\Pr(l(X_i^{G_2}) \leq t) \leq \bar{p}$.

Notice $l(X_i^{G_2})$ is approximately $\| y -\mu_1 \|_2^2$ with $y \sim N(\mu_2,\sigma^2 I)$.  By Theorem 1 in Hsu, Kakade \& Zhang(2011), for $y \sim N(\mu_2,\sigma^2 I)$,
$$\Pr\big(\| y -\mu_1 \|_2^2 \geq n\sigma^2+2\sqrt{nt}\sigma^2 + 2\sigma^2 t + \|\mu_2 - \mu_1\|_2^2(1+2\sqrt{\frac{t}{n}} ) \big) \leq \exp (-t),$$
while for the first term we can still apply $\Pr(\chi^2_n \geq t) \leq \exp(- \frac{t-\sqrt{n(2t-n)}}{2})$.
Therefore, $\Pr(l(X_i^{G_1}) \leq t)$ can be reduced to
 \begin{align*}
 \begin{split}
 \Pr( l(X_i^{G_1}) \leq t ) = & \Pr(\frac{\|x-\mu_1\|^2_2}{\hat{\sigma}_1^2} + 2 \frac{(\mu_1 - \hat{\mu}_1)'(x-\mu_1)}{{\hat{\sigma}_1^2}} + \frac{\|\mu_1 - \hat{\mu}_1\|_2^2}{\hat{\sigma}_1^2} \leq t) \\
  \geq & \Pr(\frac{\|x-\mu_1\|^2_2}{\hat{\sigma}_1^2} + 2 \frac{\|\mu_1 - \hat{\mu}_1\| \|x-\mu_1\|}{{\hat{\sigma}_1^2}} + \frac{\|\mu_1 - \hat{\mu}_1\|_2^2}{\hat{\sigma}_1^2} \leq t) \\
    = & \Pr (\frac{\|x-\mu_1\|_2}{\sigma_1} \leq -\frac{\|\mu_1 -\hat{\mu}_1\|_2}{\sigma_1} + \frac{\hat{\sigma}_1}{\sigma_1}\sqrt{t} ) \\
   \geq & 1- \exp(- \frac{(-\frac{\|\mu_1 -\hat{\mu}_1\|_2}{\sigma_1} + \frac{\hat{\sigma}_1}{\sigma_1}\sqrt{t})^2-\sqrt{n(2[-\frac{\|\mu_1 -\hat{\mu}_1\|_2}{\sigma_1} + \frac{\hat{\sigma}_1}{\sigma_1}\sqrt{t}]^2-n)}}{2}).
 \end{split}
 \end{align*}

Therefore, let $\underline{p} =  1- \exp(- \frac{(-\frac{\|\mu_1 -\hat{\mu}_1\|_2}{\sigma} + \frac{\hat{\sigma}_1}{\sigma_1}\sqrt{t})^2-\sqrt{n(2[-\frac{\|\mu_1 -\hat{\mu}_1\|_2}{\sigma} + \frac{\hat{\sigma}_1}{\sigma_1}\sqrt{t}]^2-n)}}{2})$, by the Relative Entropy Chernoff Bound,
$$ \Pr(\frac{S_{m_1}(\underline{p})}{m_1} \leq1-\beta-\frac{m_3}{m_1}) \leq \exp(-m_1 \frac{[\beta+\frac{m_3}{m_1} - \exp(- \frac{(\frac{\|\mu_1 -\hat{\mu}_1\|_2}{\sigma} + \frac{\hat{\sigma}_1}{\sigma_1}\sqrt{t})^2-\sqrt{n(2[-\frac{\|\mu_1 -\hat{\mu}_1\|_2}{\sigma} + \frac{\hat{\sigma}_1}{\sigma_1}\sqrt{t}]^2-n)}}{2})]^2 }{2 (1-\beta -\frac{m_3}{m_1})(\beta + \frac{m_3}{m_1}) } ),$$
 and suppose that
 \begin{equation}
 \label{eq18}
 \beta+\frac{m_3}{m_1} \geq  2  \exp(- \frac{(-\frac{\|\mu_1 -\hat{\mu}_1\|_2}{\sigma} + \frac{\hat{\sigma}_1}{\sigma_1}\sqrt{t})^2-\sqrt{n(2[-\frac{\|\mu_1 -\hat{\mu}_1\|_2}{\sigma} + \frac{\hat{\sigma}_1}{\sigma_1}\sqrt{t}]^2-n)}}{2}).
 \end{equation}
 To prove \eqref{eq18}, it suffices to show the following, 
   \begin{equation}
 \label{eq3}
(- \frac{\|\mu_1 -\hat{\mu}_1\|_2}{\sigma} + \frac{\hat{\sigma}_1}{\sigma_1}\sqrt{t})^2 \geq \frac{\big(\sqrt{n}+ 2\sqrt{\log(\frac{\beta}{4} + \frac{m_3}{4m_1})^{-1}}\big)^2 +n}{2}.
 \end{equation}
For $\Pr(l(X_i^{G_2}) \leq t)$, similarly, we can get
\begin{align*}
\begin{split}
\Pr(  l(X_i^{G_2}) \leq t ) & =  \Pr(\frac{\|x-\mu_1\|^2_2}{\hat{\sigma}_1^2} + 2 \frac{(\mu_1 - \hat{\mu}_1)'(x-\mu_1)}{{\hat{\sigma}_1^2}} + \frac{\|\mu_1 - \hat{\mu}_1\|_2^2}{\hat{\sigma}_1^2} \leq t) \\
& \leq  \Pr(\frac{\|x-\mu_1\|^2_2}{\hat{\sigma}_1^2} - 2 \frac{\|\mu_1 - \hat{\mu}_1\|\|x-\mu_1\|}{{\hat{\sigma}_1^2}} + \frac{\|\mu_1 - \hat{\mu}_1\|_2^2}{\hat{\sigma}_1^2} \leq t) \\
& = \Pr(\|x-\mu_1\|_2 \leq \hat{\sigma}_1\sqrt{t} + \|\mu_1 - \hat{\mu}_1\| ) \\
& = \Pr(\chi_n^2(\lambda) \leq \frac{(\hat{\sigma}_1\sqrt{t} + \|\mu_1 - \hat{\mu}_1\|)^2}{\sigma_1^2}) \\
& \leq \bar{p}.
\end{split}
\end{align*}
thus we can choose $\bar{p}$ such that  $\beta \frac{m_1}{m_2} >  \bar{p} .$
 By the Chernoff-like bound,
   $$\Pr(\frac{S_{m_2}(\bar{p})}{m_2} \geq \frac{m_1 \beta}{m_2}) \leq  \exp( -m_2 \frac{[\beta \frac{m_1}{m_2} -   \exp( - ( \frac{n+\lambda - x}{2\sqrt{n+2\lambda} } )^2 ]^2}{ \beta \frac{m_1}{m_2} ( 1-\beta \frac{m_1}{m_2})  } ), $$
when
\begin{equation}
\label{eq4}
\beta \frac{m_1}{m_2} \geq 2  \exp( - ( \frac{n+\lambda - x}{2\sqrt{n+2\lambda} } )^2.
\end{equation}
where $x = (\frac{\hat{\sigma}_1\sqrt{t} + \|\mu_1 - \hat{\mu}_1\|_2 }{\sigma_1})^2$.
Then it suffices to show that 
\begin{equation}
\label{eq6}
 (\frac{\hat{\sigma}_1\sqrt{t} + \|\mu_1 - \hat{\mu}_1\|_2 }{\sigma_1})^2 \leq n+\lambda - 2\sqrt{(n+2\lambda)\log\frac{m_2}{2 m_1\beta}}.
\end{equation}
If it holds, we have $$\Pr(  l(X^{G_1})_{m_1(1-\beta) -m_3 :m_1} \geq l(X^{G_2})_{m_1\beta :m_2}) \sim O(\exp(-m_1 \beta))$$
 Therefore, let $\beta = \frac{m_2}{m_1}\eta ,$ then
\begin{align*}
\begin{split}
& \beta +\frac{m_3}{m_1} + \frac{(m_2 -m_1 \beta)^2}{2m_1} \Pr(  l(X^{G_1})_{m_1(1-\beta) -m_3 :m_1} \geq l(X^{G_2})_{m_1\beta :m_2})  \\
= &  \frac{m_2}{m_1} \eta + \frac{m_3}{m_1}  + c \frac{m_2}{m_1} \exp ( -m_2 \eta + \log m_2),
\end{split}
\end{align*}
and as long as $\eta > \frac{\log m_2}{m_2} $, we have 
$$\frac{m_2}{m_1} \eta + \frac{m_3}{m_1}  + c \frac{m_2}{m_1} \exp ( -m_2 \eta + \log m_2)= o(\frac{m_2}{m_1}).$$
This complete the proof for the upper bound of $ \Pr( l(w) \geq l(X)_{(m_1)} | w \in G_1 ) $.

To show that such choices of $\bar{p},\underline{p}$ exist, combining \eqref{eq3} \eqref{eq6}, we have
\begin{align}
\label{eq8}
\begin{split}
 n+\lambda - 2\sqrt{(n+2\lambda)\log\frac{m_2}{2m_1\beta}} & \geq (\sqrt{b}+ 2\sqrt{c})^2 ,  \\
 b& =  \frac{\big(\sqrt{n}+ 2\sqrt{\log(\frac{\beta}{4} + \frac{m_3}{4m_1})^{-1}}\big)^2 +n}{2}, \\
 c & =  \frac{\|\mu_1 -\hat{\mu}_1\|_2^2}{ \sigma_1^2},
\end{split}
\end{align}
which yields,
\begin{align*}
\begin{split}
 \lambda \geq &  4\log{\frac{2}{\eta}} + 2\log{\frac{1}{\delta}} + 2\sqrt{n\log{\frac{1}{\delta}}} + 2\sqrt{\log\frac{1}{\eta}} \sqrt{4 \log\frac{1}{\eta} + 8 \log\frac{1}{\delta}+4\sqrt{n\log\frac{1}{\delta}} +n + 2\Delta}  +\Delta, \\
\Delta =  & 4c+\sqrt{c (2n+4\sqrt{n\log\frac{1}{\delta}} + 4\log\frac{1}{\delta})} = c'_0 w_2 \sqrt{(n+\log\frac{1}{\delta})\log n },  \\
\delta = & \frac{\beta}{4} + \frac{m_3}{4m_1}.
\end{split}
\end{align*}
One could prove that for $\eta> w_3$,  $\delta = O(\eta)$. After some simplification, we get
\begin{align}
\begin{split}
 \lambda \geq &  c'_1 \log{\frac{1}{\delta}} + c'_2\sqrt{n\log{\frac{1}{\delta} } } + c'_3 \sqrt{\log\frac{1}{\eta}( \log\frac{1}{\delta} +n )}, \\
\delta = & \frac{\beta}{4} + \frac{m_3}{4m_1},
\end{split}
\end{align}
where $c'_1= 6, c'_2 = 2, c'_3= 8+ 4c'_0 $. 
If $n=o(\log \frac{1}{\delta})$, then
\begin{align}
\begin{split}
\lambda \geq (c'_1+c'_2+c'_3) \log \frac{1}{\delta}.
\end{split}
\end{align}
Similarly  we can  find a lower bound for $\Pr(  l(w) \geq l(X)_{(m_1)} | w \in G_2)\geq 1/2$. 

Hence, when
$m \geq \frac{2 \log{\frac{1}{\eta}}} {w_2 \eta^2}$, we have $ \exp(-\frac{m (w_2 \eta)^2 } {2w_2})\leq \eta,$
which bounds the probability of Bernoulli trials,
$$\Pr(| \frac{m_2}{m} - w_2 | \geq w_2 \eta ) \leq 2\exp(-\frac{m w_2 \eta^2}{2}),$$
Similarly,
with $m \geq \frac{2 \log{\frac{1}{\eta}}} {w_3}$, we can bound the probability of Bernoulli trials,
$$\Pr(| \frac{m_3}{m} - w_3 | \geq  w_3 ) \leq 2\exp(-\frac{m w_3 }{2}).$$
Hence, by union bounds, $\Pr( w \in G_1 | l(w) \geq l(X)_{(m_1)}) \leq  32 \eta $. Also,
$$m = \Omega (\frac{\frac{n}{w_2}( \log n + \log \frac{1}{\epsilon}) \log n}{\epsilon^2} ),$$
$$m \geq \frac{2 \log \frac{1}{\epsilon + \frac{w_3}{w_2} } }{w_2 ( \epsilon + \frac{w_3}{w_2})^2}, $$
$$ m \geq \frac{2 \log \frac{1}{\epsilon + \frac{w_2}{w_3} } }{w_3}. $$
Hence, 
$$ m = \Omega ( \frac{n\frac{1}{w_2}( \log n + \log \frac{1}{\epsilon}) \log n}{\epsilon^2} +  \frac{ \frac{1}{w_2} \log (\frac{1}{\epsilon} + {\frac{w_2}{w_3} }) }{\epsilon^2 + (\frac{w_3}{w_2})^2 } +  \frac{ \log (\frac{1}{\epsilon} + {\frac{w_2}{w_3} }) }{w_3} ). $$

 \section{Proof of Proposition \ref{prop_sensitive}}
 Assume that $w'_1 > w_1$. We are interested in the following sensitivity of the filter step to misestimated $w_i$ on the input. 
\begin{align*}
\begin{split}
& \big |\Pr( w \in G_1 | l(w) \geq l(X)_{(m w'_1)})  -\Pr( w \in G_1 | l(w) \geq l(X)_{(m w_1)}) \big | \\
 = & \big | \frac{\Pr( l(w) \geq l(X)_{(mw'_1)} | w \in G_1 ) \Pr( w \in G_1) }{ \sum_{k=1}^2 \Pr(  l(w) \geq l(X)_{(m w'_1)} | w \in G_k) \Pr( w \in G_k)  +  \Pr(  l(w) \geq l(X)_{(m w'_1)} | w \in N) \Pr( w \in N) } \\
 & -  \frac{\Pr( l(w) \geq l(X)_{(mw_1)} | w \in G_1 ) \Pr( w \in G_1) }{ \sum_{k=1}^2 \Pr(  l(w) \geq l(X)_{(m w_1)} | w \in G_k) \Pr( w \in G_k)  +  \Pr(  l(w) \geq l(X)_{(m w_1)} | w \in N) \Pr( w \in N) }\big |  \\
 \leq & \big | f( \Pr( l(w) \geq l(X)_{(mw_1)} | w \in G_1 ), \Pr( l(w) \geq l(X)_{(mw_1)} | w \in G_2 )) \\
 & - f( \Pr( l(w) \geq l(X)_{(mw'_1)} | w \in G_1 ), \Pr( l(w) \geq l(X)_{(mw'_1)} | w \in G_2 ) ) \big | \\
 & + \big | g(Pr( l(w) \geq l(X)_{(mw_1)} | w \in G_1 ),\Pr( l(w) \geq l(X)_{(mw_1)} | w \in G_2 ), \Pr( l(w) \geq l(X)_{(mw_1)} | w \in N ) ) \\
 & -  g(Pr( l(w) \geq l(X)_{(mw'_1)} | w \in G_1 ),\Pr( l(w) \geq l(X)_{(mw'_1)} | w \in G_2 ), \Pr( l(w) \geq l(X)_{(mw'_1)} | w \in N ) ) \big | \\
   \leq & \max |\frac{\partial f}{\partial x}(x,y)|   | \Pr(l(w)> l(X)_{(m w'_1)}| w\in G_1) -\Pr(l(w)> l(X)_{(m w_1)}| w\in G_1) | \\
 & +  \max |\frac{\partial f}{\partial y} (x,y)|  | \Pr(l(w)> l(X)_{(m w'_1)}| w\in G_2) -\Pr(l(w)> l(X)_{(m w_1)}| w\in G_2)| \\
 & + \max |\frac{\partial g}{\partial x}(x,y,z)|   | \Pr(l(w)> l(X)_{(m w'_1)}| w\in G_1) -\Pr(l(w)> l(X)_{(m w_1)}| w\in G_1) | \\
 & +  \max |\frac{\partial g}{\partial y} (x,y,z)|  | \Pr(l(w)> l(X)_{(m w'_1)}| w\in G_2) -\Pr(l(w)> l(X)_{(m w_1)}| w\in G_2)| \\
 & + 2 \max | \frac{\partial g}{\partial z}(x,y,z) |  \\
 \leq &  \frac{2w_1 }{w_2 }  | \Pr(l(w)> l(X)_{(m w'_1)}| w\in G_1) -\Pr(l(w)> l(X)_{(m w_1)}| w\in G_1) | \\
 & +\frac{2w_1 \eta}{w_2 }  | \Pr(l(w)> l(X)_{(m w'_1)}| w\in G_2) -\Pr(l(w)> l(X)_{(m w_1)}| w\in G_2)|  \\
 & + \frac{4w_1 w_3}{w_2^2 }| \Pr(l(w)> l(X)_{(m w'_1)}| w\in G_1) -\Pr(l(w)> l(X)_{(m w_1)}| w\in G_1) | \\
 & + \frac{16 w_1 w_3 \eta }{w_2 } | \Pr(l(w)> l(X)_{(m w'_1)}| w\in G_2) -\Pr(l(w)> l(X)_{(m w_1)}| w\in G_2)|\\
 & + \frac{8w_1 w_3\eta}{w_2 } ,  
   \end{split}
 \end{align*}
  where
 \begin{align*}
 \begin{split}
f(x,y) = & \frac{w_1 x}{ w_1 x + w_2 y}. \\
 g(x,y,z) = & \frac{w_1 w_3 xz} {( w_1 x+ w_2 y)(w_1 x + w_2 y + w_3 z)} .
\end{split}
\end{align*}
The third  inequality is due to the intermediate value theorem. 

Now the problem reduces to comparing
\begin{enumerate}[(a)]
\item $\Pr(l(w)> l(X)_{(m w'_1)}| w\in G_1) $ and $\Pr(l(w)> l(X)_{(m w_1)}| w\in G_1) $ .
\item $\Pr(l(w)> l(X)_{(m w'_1)}| w\in G_2) $ and $\Pr(l(w)> l(X)_{(m w_1)}| w\in G_2) $ .
\end{enumerate}
For (a), 
\begin{align*}
\begin{split}
& \Pr( l(X)_{(m w'_1)}  > l(w)>  l(X)_{(m w_1)} | w\in G_1)  \\
\leq & \Pr( l(X)_{(m w'_1)}  > l(w)>  l(X)_{(m w_1)}, l(w) \leq l(X^{G_1})_{(m_1(1-\tau))} | w\in G_1) + \tau \\
\leq & \Pr( l(w)> l(X^{G_1\cup G_2})_{(mw_1-m_3)} , l(w) \leq l(X^{G_1})_{(m_1(1-\tau))}| w\in G_1) \\
&   - \Pr(l(w) > l(X^{G_1\cup G_2})_{(mw'_1)} , l(w) \leq l(X^{G_1})_{(m_1(1-\tau))}|w\in G_1)   \\
\leq & \sum_{n_1 + n_2=mw_1 - m_3} \Pr( l(w) > l(X^{G_1})_{(n_1) }, l(w) > l(X^{G_2})_{(n_2)}, l(w) \leq l(X^{G_1})_{(m_1(1-\tau))} | w\in G_1) \\
&  -\sum_{n'_1 + n'_2=mw'_1 } \Pr( l(w) > l(X^{G_1})_{(n'_1) },l(w) \leq l(X^{G_1})_{(m_1(1-\tau))}, l(w) > l(X^{G_2})_{(n'_2)} | w\in G_1)  + 2\tau \\
\leq &  \sum_{n_1 + n_2=mw_1 - m_3}\sum_{j= n_1}^{m_1(1-\tau)} \frac{1}{m_1} \Pr( l(X^{G_1})_{(j)} > l(X^{G_2})_{(n_2)})  - \sum_{n'_1 + n'_2=mw'_1}\sum_{q = n'_1}^{m_1(1-\tau)} \frac{1}{m_1} \Pr( l(X^{G_1})_{(q) } > l(X^{G_2})_{(n'_2)})  +2 \tau \\
\leq & \sum_{n_1 = mw_1 -m_3 -m_2}^{mw_1-m_3} \sum_{j=n_1}^{m_1(1-\tau)}  \frac{1}{m_1} \Pr( l(X^{G_1})_{(j) } > l(X^{G_2})_{(mw_1 - m_3 - n_1)}) - \sum_{n'_1 = mw'_1 - m_2}^{mw_1} \sum_{j=n'_1}^{m_1(1-\tau)}  \frac{1}{m_1} \Pr( l(X^{G_1})_{(j) } > l(X^{G_2})_{(mw'_1-n'_1 )})  + 2 \tau. 
\end{split}
\end{align*}
Let $n'_1 = l_1 + m(w'_1 - w_1) +m_3, $ then it is equivalent to 
\begin{align*}
\begin{split}
&  \sum_{n_1 = mw_1 -m_3 -m_2}^{(mw_1-m_3)} \sum_{j=n_1}^{m_1(1-\tau)} \frac{1}{m_1}  \Pr( l(X^{G_1})_{(j) } > l(X^{G_2})_{(mw_1 - m_3 - n_1)}) \\
& -  \sum_{l_1  = mw_1 - m_2-m_3}^{mw_1 - m(w'_1 - w_1) -m_3 } \sum_{j=l_1 + m(w'_1 - w_1) +m_3}^{m_1(1-\tau)}  \frac{1}{m_1}  \Pr( l(X^{G_1})_{(j) } > l(X^{G_2})_{(mw_1 -m_3 -l_1 )}) + 2 \tau\\
= & \sum_{i=mw_1 - m_3 -m_1(1-\tau)}^{m(w'_1 - w_1)} \sum_{j=mw_1 - m_3-i}^{m_1(1-\tau) } \frac{1}{m_1} \Pr( l(X^{G_1})_{(j) } > l(X^{G_2})_{(i)})  +  \sum_{i=mw'_1 - m_1(1-\tau)}^{m_2} \sum_{j=mw_1 -m_3 -i}^{mw'_1 -i} \frac{1}{m_1}  \Pr( l(X^{G_1})_{(j) } > l(X^{G_2})_{(i)}) + 2\tau \\
\leq & \frac{[m_1(1-\tau) -mw_1 + m_3 + m(w'_1 - w_1)]^2 }{2m_1}\Pr(l(X^{G_1})_{(m_1(1-\tau))} > l(X^{G_2})_{mw_1 -m_3 -m_1(1-\tau) }) \\
& + \frac{[m(w'_1 - w_1)+m_3][m_2 +m_1(1-\tau) - m w'_1]}{m_1} \Pr(  l(X^{G_1})_{(m_1(1-\tau))}> l(X^{G_2})_{(mw'_1-m_1(1-\tau))}) +2 \tau \\
\leq c_5 & \frac{[m(w'_1 - w_1) + m_3 - m_1 \tau ]^2}{2m_1}  \Pr(l(X^{G_1})_{(m_1(1-\tau))} > l(X^{G_2})_{(m_1 \tau -m_3  )}) \\
& + \frac{[m(w'_1 - w_1) + m_3][m_2 - m_1 \tau - m(w'_1 - w_1)] }{m_1}  \Pr(  l(X^{G_1})_{(m_1(1-\tau))}> l(X^{G_2})_{(m(w'_1-w_1)+ m_1\tau)}  ) + 2\tau  \\
= & \frac{[m(w'_1 - w_1) + m_3 - m_1 \tau ]^2}{2m_1}c_5 \exp(-m_1\tau) + \frac{[m(w'_1 - w_1) + m_3][m_2 - m_1 \tau - m(w'_1 - w_1)] }{m_1} c_5  \exp(-m_1\tau) + 2\tau,
\end{split}
\end{align*}
under Lemma \ref{lemma1s}. Let $\tau= \min \{  \frac{w_3}{w_1}+ \frac{w'_1 - w_1}{ w_1} , \frac{w_2}{w_1} - \frac{w'_1 - w_1}{ w_1}\}$, we have 
\begin{align*}
\begin{split}
 & \Pr( l(X)_{(m w'_1)}  > l(w)>  l(X)_{(m w_1)} | w\in G_1) \\
= & \Pr(l(w)> l(X)_{(m w_1)}| w\in G_1)- \Pr(l(w)> l(X)_{(m w'_1)}| w\in G_1)  \\
\leq &  c_5 m_1\tau \exp(-m_1 \tau) (\frac{w_2}{w_1} - \frac{w'_1 - w_1}{w_1} -\tau )  +c_5 m_1 \exp(-m_1 \tau) (\frac{w_3}{w_1}+ \frac{w'_1 - w_1}{ w_1}  - \tau)^2 + 2 \tau \\
\leq  & c_5 m_1 \exp(-m_1 \tau)  \frac{w_2+w_3}{w_1} ( \frac{w_3}{w_1}+ \frac{w'_1 - w_1}{ w_1} ) + 2(\frac{w_3}{w_1}+ \frac{w'_1 - w_1}{ w_1} )\\
\leq & c_6 (\frac{w_3}{w_1}+ \frac{w'_1 - w_1}{ w_1} ),
\end{split}
\end{align*}
for some $c_6>0$.

Similarly, the following argues for the case  when $w\in G_2$.
\begin{align*}
\begin{split}
 & \Pr( l(X)_{(m w'_1)}  > l(w)>  l(X)_{(m w_1)} | w\in G_2) \\
= & \Pr( l(X)_{ mw'_1 } > l(w),  | w\in G_2) - \Pr(l(X)_{(mw_1)} |w \in G_2) \\
\leq & \Pr(l(X^{G_1 \cup G_2} )_{(mw'_1 )} > l(w), l(w) \geq l(X^{G_2})_ {(mw'_1 -m_1 (1-\tau)) } |w \in G_2) \\
& - \Pr(l(X^{G_1 \cup G_2} )_{(mw_1-m_3)} > l(w), l(w) \geq l(X^{G_2})_ {(m_1 \tau + mw_1 - m_1 -m_3) } |w \in G_2) \\
& + \frac{mw'_1 -m_1 (1-\tau)}{m_2 } + \frac{m_1 \tau + mw_1 - m_1 -m_3}{m_2} \\
\leq & \sum_{n'_1 + n'_2 = m-m_3-mw'_1, n'_1 \geq m_1 \tau } \Pr(l(w)< l(X^{G_1})_{(m_1 - n'_1) } , l(w) < l(X^{G_2})_ {(m_2- n'_2 )}, l(w) \geq l(X^{G_2})_ {(mw'_1 -m_1 (1-\tau)) } | w\in G_2) \\
 & -  \sum_{n_1 + n_2 = m-mw_1, n_1 \geq m_1 \tau } \Pr(l(w)< l(X^{G_1})_{(m_1 - n_1) } , l(w) < l(X^{G_2})_ {(m_2- n_2 )}, l(w) \geq l(X^{G_2})_ {(m_1 \tau + mw_1 - m_1 -m_3) } | w\in G_2)\\
 &  +  \frac{mw'_1 -m_1(1-\tau)}{m_2}  + \frac{mw_1 - m_1(1-\tau) - m_3}{m_2} \\
 = & \sum_{n'_1 = m_1 \tau}^{m- m_3 -mw'_1} \sum_{j=0}^{ mw'_1 + n'_1 -m_1 } \frac{1}{m_2} \Pr(l(X^{G_2})_{(j)} < l(X^{G_1})_{(m_1 - n'_1 )})  - \sum_{n_1= m_1\tau}^{m-mw_1} \sum_{j=0}^{mw_1-m_1 -m_3 + n_1)} \frac{1}{m_2} \Pr(l(X^{G_2})_{(j)} < l(X^{G_1})_{(m_1 - n_1) } ) \\
 &  +  \frac{mw'_1 -m_1(1-\tau)}{m_2}  + \frac{mw_1 - m_1(1-\tau) - m_3}{m_2} \\
  = & \sum_{i=mw'_1 -m_2}^{m(1-\tau)} \sum_{j=0}^{mw'_1 - i} \frac{1}{m_2} \Pr( l(X^{G_1})_{(i)} > l(X^{G_2})_{(j)} )  - \sum_{i= m_1 -m+mw_1}^{m_1(1-\tau)} \sum_{j=0}^{mw_1 -m_3 - i } \frac{1}{m_2} \Pr( l(X^{G_1})_{(i)} > l(X^{G_2})_{(j)} )   \\
  & +  \frac{mw'_1 -m_1(1-\tau)}{m_2}  + \frac{mw_1 - m_1(1-\tau) - m_3}{m_2}\\
  =& \sum_{i=mw'_1 -m_2}^{m_1(1-\tau)} (\sum_{j=0}^{mw'_1 - i} - \sum_{j=0}^{mw_1 - m_3 -i} ) \frac{1}{m_2} \Pr( l(X^{G_1})_{(i)} > l(X^{G_2})_{(j)} )    - \sum_{i=m_1 - m + mw_1}^{mw'_1 - m_2} \sum_{j=0}^{mw_1 - m_3 -i} \frac{1}{m_2} \Pr( l(X^{G_1})_{(i)} > l(X^{G_2})_{(j)} ) \\
  & +  \frac{mw'_1 -m_1(1-\tau)}{m_2}  + \frac{mw_1 - m_1(1-\tau) - m_3}{m_2}\\
  \leq & \sum_{i=mw'_1 -m_2}^{m_1(1-\tau)} \sum_{j=mw_1 -m_3 -i}^{mw'_1 -i} \frac{1}{m_2} \Pr( l(X^{G_1})_{(i)} > l(X^{G_2})_{(j)} ) + \frac{mw'_1 -m_1(1-\tau)}{m_2}  + \frac{mw_1 - m_1(1-\tau) - m_3}{m_2} \\
  = & \sum_{i=m_1 \tau}^{m_2 + m_1 -mw'_1} \sum_{j = mw_1 -m_3 -m_1 +i}^{mw'_1 -m_1 +i}  \frac{1}{m_2} \Pr( l(X^{G_1})_{(m_1 - i)} > l(X^{G_2})_{(j)} ) +  \frac{mw'_1 -m_1(1-\tau)}{m_2}  + \frac{mw_1 - m_1(1-\tau) - m_3}{m_2} \\
 =  &  \sum_{i=m_1 \tau}^{m_2 + m_1 -mw'_1} \sum_{j = mw_1 -m_3 -m_1 +i}^{mw'_1 -m_1 +i}  \frac{1}{m_2} c_1 \exp( -m_1 \tau ) +  \frac{mw'_1 -m_1(1-\tau)}{m_2}+ \frac{mw_1 - m_1(1-\tau) - m_3}{m_2} \\
 = & \frac{(m_1(1-\tau) + m_2 - mw'_1)(m(w'_1 - w_1) + m_3)}{m_2}  c_1 \exp(-m_1 \tau)  +  \frac{mw'_1 -m_1(1-\tau)}{m_2}+ \frac{mw_1 - m_1(1-\tau) - m_3}{m_2} \\
 = & \frac{w'_1 - w_1+ w_3}{w_2} ( \frac{w_2}{w_1} -\tau - \frac{w'_1-w_1}{w_1}) m_1 c_1 \exp(-m_1 \tau) + \frac{w'_1 - w_1 - w_3}{w_2} + 2 \frac{w_1}{w_2} \tau \\
 = &  3\frac{w'_1 - w_1 }{w_1} + c_1m_1 \exp(-m_1 \tau)  \frac{w'_1 - w_1 + w_3}{w_1},
\end{split}
\end{align*}
where $\tau=\frac{w'_1 - w_1 (1-\frac{w_2}{w_1}) -w_3}{2 w_1}$.
Hence, we can find $c_7, c_8 >0$ such that, 
 $$ \big | \Pr(l(w)> l(X)_{(m w'_1)}| w\in G_1)- \Pr(l(w)> l(X)_{(m w_1)}| w\in G_1)\big | \leq c_7 \frac{w'_1 - w_1 + w_3}{w_1},  $$ 
 $$\big | \Pr(l(w)> l(X)_{(m w'_1)}| w\in G_2) - \Pr(l(w)> l(X)_{(m w_1)}| w\in G_2) \big | \leq  c_8 \frac{w'_1 - w_1 + w_3}{w_1}.  $$
 
 Therefore,
 \begin{align*}
\begin{split}
&  | \Pr( w \in G_1 | l(w) \geq l(X)_{(m w'_1)})  -  \Pr( w \in G_1 | l(w) \geq l(X)_{(m w_1)}) |  \\
 \leq &  \frac{2w_1 }{w_2 }  | \Pr(l(w)> l(X)_{(m w'_1)}| w\in G_1) -\Pr(l(w)> l(X)_{(m w_1)}| w\in G_1) | \\
 & +\frac{2w_1 \eta}{w_2 }  | \Pr(l(w)> l(X)_{(m w'_1)}| w\in G_2) -\Pr(l(w)> l(X)_{(m w_1)}| w\in G_2)|  \\
 & + \frac{4w_1 w_3}{w_2^2 }| \Pr(l(w)> l(X)_{(m w'_1)}| w\in G_1) -\Pr(l(w)> l(X)_{(m w_1)}| w\in G_1) | \\
 & + \frac{16 w_1 w_3 \eta }{w_2 } | \Pr(l(w)> l(X)_{(m w'_1)}| w\in G_2) -\Pr(l(w)> l(X)_{(m w_1)}| w\in G_2)|\\
 & + \frac{8w_1 w_3\eta}{w_2 } \\
 \leq &  (\frac{2w_1 }{w_2 }+\frac{4w_1 w_3}{w_2^2})  | \Pr(l(w)> l(X)_{(m w'_1)}| w\in G_1) -\Pr(l(w)> l(X)_{(m w_1)}| w\in G_1) | \\
 & + (\frac{2w_1\eta }{w_2 }+\frac{ 16w_1 w_3\eta}{w_2}) | \Pr(l(w)> l(X)_{(m w'_1)}| w\in G_2) -\Pr(l(w)> l(X)_{(m w_1)}| w\in G_2)| + \frac{8w_1 w_3\eta}{w_2}  \\
 \leq & c_9 \frac{w'_1 - w_1 + w_3}{w_2},
  \end{split}
 \end{align*}
 for some $c_9 >0$.
 
 Hence, when reporting $w'_1$, the malicious noise of weight accounts for at most,  
 $$ \eta' = \eta + c_9  \frac{w'_1 - w_1 + w_3}{w_2}, $$
and after the agnostic recovery of the second mean, $\| \hat{\mu}_2(w'_1) - \mu_2\|$ is bounded by 
\begin{align*}
\begin{split}
O(\eta + c_9  \frac{w'_1 - w_1 + w_3}{w_2}  + \frac{w_3}{1-w'_1} +  \epsilon ) \sigma \sqrt{\log n } ),
\end{split}
\end{align*}
given the same condition of sampling complexity regarding $\eta$. Here $\frac{w_3}{1-w'_1} $ is the upper bound for the weight of malicious noise in the new input in the second run.
In the spherical case, $\| \hat{\mu}_2(w_1) - {\mu}_2\|$ is bounded by $O(\eta +  \frac{ w_3}{w_2}  +  \epsilon ) \sigma \sqrt{\log n } )$.
Hence, 
$\| \hat{\mu}_2(w'_1) - \hat{\mu}_2(w_1)\|  $ is bounded by
\begin{align*}
\begin{split}
O( \frac{w'_1 - w_1 + w_3}{w_2}  + \frac{w_3}{1-w'_1} +  \epsilon ) \sigma \sqrt{\log n }).
\end{split}
\end{align*}
 Similarly, one could prove  the same for the case $w'_1 < w_1$. Therefore, 
 $$\| \hat{\mu}_2(w'_1) - \hat{\mu}_2(w_1)\|  = O( \frac{|w'_1 - w_1| + w_3}{w_2}  + \frac{w_3}{1-w'_1} +  \epsilon ) \sigma \sqrt{\log n }). $$
 
 \end{document}